\documentclass[12pt]{report} 
\usepackage{latexsym}
\newcommand{\hs}{\hspace{0.3cm}}
\newcommand{\insp}{\hspace*{1cm}}

\newcommand{\ds}{\displaystyle}

\usepackage{amsfonts}

\usepackage{epsfig}

\newcommand{\reg}{\mbox{\scriptsize{reg}}}

\newcommand{\codim}{\mbox{codim}}
\newcommand{\trdeg}{\mbox{trdeg}}
\newcommand{\sing}{\mbox{\scriptsize{sing}}}

\newcommand{\rank}{\mbox{rank}}

\newcommand{\IZ}{Z\!\!\!Z}

\newcommand{\IC}{I\!\!\!\!C}

\newcommand{{\CB}}{\cal B}

\newcommand{{\CC}}{\cal C}

 1
 1

\def\picture#1 by #2 (#3){
\vbox to #2{
\hrule width #1 height 0pt depth 0pt
\vfill
\special{picture #3}}}

\begin{document}
\begin{center}
{\bf \Large The Poisson center and polynomial, maximal Poisson commutative subalgebras, especially for nilpotent Lie algebras of dimension at most seven}\\
\ \\
{\bf \large Alfons I. Ooms}\\
\ \\
{\it Mathematics Department, Hasselt University, Agoralaan, Campus Diepenbeek, 3590 Diepenbeek, Belgium\\
E-mail address: alfons.ooms@uhasselt.be}
\end{center}

To the memory of Andr\'e Cerezo\footnote{http://math.unice.fr/$\sim$frou/AC.html} (1945-2003)\\
\ \\
{\bf Key words:} Poisson center, Poisson commutative subalgebra, nilpotent Lie algebra, enveloping algebra.\\
AMS classification: 17B35.\\
\ \\
{\bf Abstract.} Let $\mathfrak{g}$ be a finite dimensional Lie algebra over an algebraically closed field $k$ of characteristic zero.  We collect some general results on the Poisson center of $S(\mathfrak{g})$, including some simple criteria regarding its polynomiality, and also on certain Poisson commutative subalgebras of $S(\mathfrak{g})$.  These facts are then used to finish the study of [05,5], i.e. to give an explicit description for the Poisson center of all indecomposable, nilpotent Lie algebras of dimension at most seven.  Among other things, we also provide a polynomial, maximal Poisson commutative subalgebra of $S(\mathfrak{g})$, enjoying additional properties.  As a by-product we show that a conjecture by Milovanov is valid in this situation.  These results easily carry over to the enveloping algebra $U(\mathfrak{g})$.\\
\ \\
{\bf \large 1. Introduction}\\
{\bf 1.1} Let $k$ be an algebraically closed field of characteristic zero and let $\mathfrak{g}$ be a Lie algebra over $k$ with basis $x_1, \ldots, x_n$. For each $\xi \in \mathfrak{g}^\ast$ we consider its stabilizer
\begin{eqnarray*}
\mathfrak{g}(\xi) = \{x \in \mathfrak{g}\mid \xi ([x,y]) = 0\ \mbox{for all}\ y \in \mathfrak{g}\}
\end{eqnarray*}
The minimal value of $\dim \mathfrak{g}(\xi)$ is called the index of $\mathfrak{g}$ and is denoted by $i(\mathfrak{g})$ [D5, 1.11.6; TY, 19.7.3].  Put $c(\mathfrak{g}) = (\dim \mathfrak{g} + i(\mathfrak{g}))/2$.  This integer will play an important role throughout this paper.  An element $\xi \in \mathfrak{g}^\ast$ is called regular if $\dim \mathfrak{g}(\xi) = i(\mathfrak{g})$.  The set $\mathfrak{g}^\ast_{\reg}$ of all regular elements of $\mathfrak{g}^\ast$ is an open dense subset of $\mathfrak{g}^\ast$.\\
We put $\mathfrak{g}^\ast_{\sing} = \mathfrak{g}^\ast \backslash \mathfrak{g}^\ast_{\reg}$.  Clearly, $\codim\ \mathfrak{g}^\ast_{\sing} \geq 1$.  Following [JS] we call $\mathfrak{g}$ singular if equality holds and nonsingular otherwise.  For instance, any semi-simple Lie algebra $\mathfrak{g}$ is nonsingular since $\codim\ \mathfrak{g}^\ast_{\sing} = 3$.\\
\ \\
{\bf 1.2.} \hs We equip the symmetric algebra $S(\mathfrak{g})$ with its natural structure of Poisson algebra.  Its center is $Y(\mathfrak{g}) = S(\mathfrak{g})^{\mathfrak{g}}$, the subalgebra of the invariant polynomials of $S(\mathfrak{g})$.  For brevity we call $\mathfrak{g}$ coregular if $Y(\mathfrak{g})$ is a polynomial algebra.  In recent years considerable attention has been devoted to when this occurs and also to the existence of certain Poisson commutative subalgebras of $S(\mathfrak{g})$ [Bo2, FJ1, FJ2, J3, J4, J5, J6, JL, JS, PY, Sa, T, Y].  In this paper we will give a few simple criteria for coregularity (Theorems 13 and 26, Examples 27 and 28). These will be applied frequently in section 3.  Of the same nature, but more profound, is the sum rule (Theorem 9 and Remark 10), a formula involving the sum of the degrees of the homogeneous generators of $Y(\mathfrak{g})$.  Joseph and Shafrir obtained in [JS, 5.7] the following converse, which is the key ingredient used in part I of section 3.  It extends a result of Panyushev, Premet and Yakimova [PPY, Theorem 1.2] (see also [P, Theorem 1.2]) who treated the nonsingular, algebraic case.\\
\ \\
{\bf Theorem 14.} Assume $\trdeg_k Y(\mathfrak{g}) = i(\mathfrak{g})$. Let $f_1,\ldots, f_r \in Y(\mathfrak{g})$, $r = i(\mathfrak{g})$, be homogeneous, algebraically independent polynomials  such that
$$\sum\limits_{i=1}^r \deg f_i \leq c(\mathfrak{g}) - \deg p_{\mathfrak{g}}$$
where $p_{\mathfrak{g}}$ is the fundamental semi-invariant of $\mathfrak{g}$ (see definition 4).  Then $Y({\mathfrak{g}})$ is freely generated by $f_1,\ldots, f_r$.\\
\ \\
Next, $Y({\mathfrak{g}}) \subset S(F(\mathfrak{g}))$ (Remark 8), where $F(\mathfrak{g}) = \sum\limits_{\xi \in \mathfrak{g}_{\reg}^\ast} \mathfrak{g}(\xi)$ is the Frobenius semi-radical of $\mathfrak{g}$ (2.5), which also plays a useful role in this paper.\\
\ \\
{\bf 1.3} Let $A$ be a Poisson commutative subalgebra of $S(\mathfrak{g})$.  Then it is well-known that $\trdeg_k(A) \leq c(\mathfrak{g})$.  $A$ is called complete if equality holds and strongly complete if it is also a maximal Poisson commutative subalgebra.  According to Sadetov there always exists a complete Poisson commutative subalgebra of $S(\mathfrak{g})$ [Sa].  For example, suppose $\mathfrak{g}$ admits a commutative Lie subalgebra $\mathfrak{h}$ of $\mathfrak{g}$ with $\dim \mathfrak{h} = c(\mathfrak{g})$, i.e. $\mathfrak{h}$ is a commutative polarization (CP) of $\mathfrak{g}$.  Then $S(\mathfrak{h})$ is a polynomial, strongly complete subalgebra of $S(\mathfrak{g})$ and its quotient field $R(\mathfrak{h})$ is a maximal Poisson commutative subfield of $R(\mathfrak{g})$.  If in addition $\mathfrak{h}$ is also an ideal of $\mathfrak{g}$ then $\{\mathfrak{g}, S(\mathfrak{h})\} \subset S(\mathfrak{h})$ (1 of Examples 19).\\
\ \\
In the nilpotent case we follow an approach by Mich\`ele Vergne [V], namely we consider an increasing sequence of ideals $\mathfrak{g}_i$ of $\mathfrak{g}$, such that $\dim \mathfrak{g}_i = i$, $i : 0,\ldots, n$ and we let $V(\mathfrak{g})$ be the subalgebra of $S(\mathfrak{g})$ generated by the union of all $Y(\mathfrak{g}_i)$.\\
Then $V(\mathfrak{g})$ is complete and its Poisson commutant $M = V(\mathfrak{g})'$ is a strongly complete Poisson commutative subalgebra of $S(\mathfrak{g})$.  Also, its quotient field $Q(M)$ is a maximal Poisson commutative subfield of $R(\mathfrak{g})$ and $\{\mathfrak{g},M\} \subset M$ (2 of Examples 19).\\
\ \\
A general criterion by Bolsinov asserts that for any $\mathfrak{g}$ and $\xi \in \mathfrak{g}_{\reg}^\ast$ the Mishchenko-Fomenko subalgebra $Y_\xi(\mathfrak{g})$ (see 2.8 for its definition) is complete if and only if $\mathfrak{g}$ is nonsingular [Bo1, Bo2].  In 2.8 we give a counterexample to this and also state the improved version (Theorem 21) by Joseph and Shafrir [JS, 7.2].\\
\ \\
{\bf 1.4.} If $\mathfrak{g}$ admits a CP $\mathfrak{h}$ then $F(\mathfrak{g})$ is commutative (since $F(\mathfrak{g}) \subset \mathfrak{h}$ [O3,  p.710]).  We will examine under which circumstances the converse holds.  Assume that $\trdeg_k Y(\mathfrak{g}) = i(\mathfrak{g})$ and that $\mathfrak{g}$ is nonsingular.  Then $c(F(\mathfrak{g})) = c(\mathfrak{g})$.  In particular, any complete Poisson commutative subalgebra of $S(F(\mathfrak{g}))$ is then also complete in $S(\mathfrak{g})$.  Moreover, if $F(\mathfrak{g})$ is commutative, then it is the only CP of $\mathfrak{g}$ (2 of Theorem 22).\\
See also Remark 23 and 3.4.4.\\
\ \\
{\bf 1.5.} In section 3 we will finish the study of [O5,5] $(k = \IC)$.  Namely, it remains to give the explicit description of the Poissoçn center $Y(\mathfrak{g})$ for each of the 77 7-dimensional indecomposable nilpotent Lie algebras $\mathfrak{g}$ for which $i(\mathfrak{g}) > \rank\ \mathfrak{g}$ (the latter being the dimension of a maximal torus inside $\mbox{Der} \ \mathfrak{g}$).  However, we will not list the 23 minimal generators of $Y(\mathfrak{g})$, where $\mathfrak{g}$ is the 7-dimensional standard filiform Lie algebra, because of the enormous size of most of these generators.  They can be found in the unpublished manuscript by Andr\'e Cerezo [Ce4], where their description takes about 20 pages (see Example 27).  For each Lie algebra of the list we will also exhibit a polynomial, strongly complete Poisson commutative subalgebra $M$ of $S(\mathfrak{g})$ such that $\{\mathfrak{g},M\} \subset M$ and such that $Q(M)$ is a maximal Poisson commutative subfield of $R(\mathfrak{g})$, by using the methods of Examples 19.  Moreover, the generators of $M$ have degrees at most 2.  This implies that a conjecture by Milovanov (2.9) holds for all nilpotent Lie algebras of dimension at most 7.  Finally, we will explain why these results can be interpreted directly as results for the enveloping algebra $U(\mathfrak{g})$, by using Proposition 18 and Examples 19.\\
\ \\
{\bf 2. Preliminaries and general results}\\
Let $k$ be an algebraically closed field of characteristic zero and let $\mathfrak{g}$ be a Lie algebra over $k$ with basis $x_1, \ldots, x_n$.\\
\ \\
{\bf 2.1 The Poisson algebra $S(\mathfrak{g})$ and its center}\\
The symmetric algebra $S(\mathfrak{g})$, which we identify with $k[x_1,\ldots, x_n]$, has a natural Poisson algebra structure, the Poisson bracket of $f, g \in S(\mathfrak{g})$ given by:
$$\{f,g\} = \sum\limits_{i=1}^n \sum\limits_{j=1}^n [x_i,x_j] \ds\frac{\partial f}{\partial x_i} \ds\frac{\partial g}{\partial x_j}$$
In particular, $S(\mathfrak{g}),\{,\}$ is a Lie algebra for which $\mathfrak{g}$ is a Lie subalgebra since for any two elements $x, y \in \mathfrak{g}$ we have that $\{x,y\} = [x,y]$.  Also, for all $f,g,h \in S(\mathfrak{g})$:
$$\{f, gh\} = \{f,g\}h + g\{f,h\} \insp (\ast)$$
It now easily follows that the center of $S(\mathfrak{g}),\{,\}$ is equal to
$$\{f \in S(\mathfrak{g})\mid \{x,f\} = 0\ \ \forall x \in \mathfrak{g}\}$$
and since $\{x,f\} = \mbox{ad}\ x(f)$ this clearly coincides with $Y(\mathfrak{g}) = S(\mathfrak{g})^{\mathfrak{g}}$, the subalgebra of invariant polynomials of $S(\mathfrak{g})$.  Furthermore, a subalgebra $A$ of $S(\mathfrak{g})$ is said to be Poisson commutative if $\{f,g\} = 0$ for all $f,g \in A$.\\
Finally, the Poisson bracket has a unique extension to the quotient field $R(\mathfrak{g})$ of $S(\mathfrak{g})$ such that $(\ast)$ holds in $R(\mathfrak{g})$.  It follows that $R(\mathfrak{g}),\{,\}$ is a Lie algebra with center $R(\mathfrak{g})^{\mathfrak{g}}$, the subfield of rational invariants of $R(\mathfrak{g})$. $R(\mathfrak{g})$ is called the rational Poisson algebra [V, p. 311].\\
\ \\
{\bf 2.2 Two formulas involving the index of $\mathfrak{g}$}\\
First, we recall from [D5, 1.14.13] that
$$i(\mathfrak{g}) = \dim \mathfrak{g} - \rank_{R(\mathfrak{g})} ([x_i, x_j])$$
Next, denote by $Z(D(\mathfrak{g}))$ the center of the quotient division ring $D(\mathfrak{g})$ of the enveloping algebra $U(\mathfrak{g})$ of $\mathfrak{g}$.\\
\ \\
{\bf Theorem 1}
$$\trdeg_k(R(\mathfrak{g})^{\mathfrak{g}}) = \trdeg_k(Z(D(\mathfrak{g}))) \leq i(\mathfrak{g})$$
Moreover, equality occurs if one of the following conditions is satisfied:
\begin{itemize}
\item[(1)] $\mathfrak{g}$ is algebraic [RV,O1].
\item[(2)] $\mathfrak{g}$ has no proper semi-invariants (in $S(\mathfrak{g})$ or equivalently in $U(\mathfrak{g})$) [OV, Proposition 4.1].
\end{itemize}

{\bf 2.3 The fundamental semi-invariant $p_{\mathfrak{g}}$}\\
Let $\lambda \in \mathfrak{g}^\ast$.  We denote by $S(\mathfrak{g})_{\lambda}$ the set of all $f \in S(\mathfrak{g})$ such that $\mbox{ad}\ x(f) = \lambda(x)f$ for all $x \in \mathfrak{g}$.  Any element $f \in S(\mathfrak{g})_{\lambda}$ is said to be a semi-invariant w.r.t. the weight $\lambda$.  We call $f$ a proper semi-invariant if $\lambda \neq 0$.  Clearly, $S(\mathfrak{g})_{\lambda} S(\mathfrak{g})_{\mu} \subset S(\mathfrak{g})_{\lambda + \mu}$ for all $\lambda, \mu \in \mathfrak{g}^\ast$.  Let $f, g \in S(\mathfrak{g})$.  If $fg$ is a nonzero semi-invariant of $S(\mathfrak{g})$, then so are $f$ and $g$.\\
A useful link with $U(\mathfrak{g})$ is the symmetrization map, i.e. the canonical linear isomorphism $s$ of $S(\mathfrak{g})$ onto $U(\mathfrak{g})$, which maps each product $y_1 \ldots y_m$, $y_i \in \mathfrak{g}$, into $(1/m!) \sum_p y_{p(1)} \ldots y_{p(m)}$, where $p$ ranges over all permutations of $\{1,\ldots, m\}$.\\
\ \\
{\bf Remark 2.} Suppose $y_1,\ldots, y_m \in \mathfrak{g}$ commute with each other.  Then clearly,\\
$s(y_1\ldots y_m) = y_1 \ldots y_m$.\\
\ \\
$s$ is known to commute with derivations of $\mathfrak{g}$ and hence maps $S(\mathfrak{g})_{\lambda}$ onto $U(\mathfrak{g})_{\lambda}$.  Taking $\lambda = 0$, the restriction
$$s : Y(\mathfrak{g}) \rightarrow Z(U(\mathfrak{g}))$$
is an algebra isomorphism if $\mathfrak{g}$ is nilpotent [D5, 4.8.12].\\
Next, take $h \in R(\mathfrak{g})$.  Then, $h \in R(\mathfrak{g})^{\mathfrak{g}}$ if and only if $h$ can be written as a quotient of two semi-invariants of $S(\mathfrak{g})$ with the same weight.  Similar properties hold in $U(\mathfrak{g})$ and $D(\mathfrak{g})$ [RV, p. 401; DNO, p. 329].\\
\ \\
{\bf Remark 3.} Assume that $\mathfrak{g}$ has no proper semi-invariants (as it is if the radical of $\mathfrak{g}$ is nilpotent).  Then $R(\mathfrak{g})^{\mathfrak{g}}$ is the quotient field of $S(\mathfrak{g})^{\mathfrak{g}} = Y(\mathfrak{g})$.  In particular,
$$\trdeg_k Y(\mathfrak{g}) = \trdeg_k R(\mathfrak{g})^{\mathfrak{g}} = i(\mathfrak{g})$$
by Theorem 1.  Also, $\mathfrak{g}$ is unimodular (i.e. $\mbox{tr}(\mbox{ad}\ x) = 0$ for all $x \in \mathfrak{g}$) by [DDV, Thm. 1.11] and its proof.\\
\ \\
{\bf Definition 4.} Put $t = \dim \mathfrak{g} - i(\mathfrak{g})$, which is the rank of the structure matrix $B = ([x_i,x_j]) \in M_n(R(\mathfrak{g}))$, where $x_1, \ldots, x_n$ is an arbitrary basis of $\mathfrak{g}$.  Assume first that $\mathfrak{g}$ is nonabelian.  Then the greatest common divisor $q_{\mathfrak{g}}$ of the $t \times t$ minors in $B$ is a nonzero semi-invariant of $S(\mathfrak{g})$
 [DNO, pp. 336-337].  If $\mathfrak{g}$ is abelian we put $q_{\mathfrak{g}} = 1$.  Next, let $p_{\mathfrak{g}}$ be the greatest common divisor of the Pfaffians of the principal $t \times t$ minors in $B$.  In particular, $\deg p_{\mathfrak{g}} \leq (\dim \mathfrak{g} - i(\mathfrak{g}))/2$.  By [OV, Lemma 2.1] $p_{\mathfrak{g}}^2 = q_{\mathfrak{g}}$ up to a nonzero scalar multiplier.  We call $p_{\mathfrak{g}}$ the fundamental semi-invariant of $S(\mathfrak{g})$ (instead of $q_{\mathfrak{g}}$ as we did in [OV, p. 309]).  This corresponds with Dixmier's notion of fundamental semi-invariant of $U(\mathfrak{g})$ where $\mathfrak{g} = af(n,\IC)$ [D4; DNO, p. 345].\\
\ \\
{\bf Remark 5.} [OV, p. 307]
\begin{center}
$\mathfrak{g}$ is singular if and only if $p_{\mathfrak{g}} \notin k$
\end{center}

{\bf Example 6.}  Let $\mathfrak{g}$ be a nonabelian Lie algebra with center $Z(\mathfrak{g})$.  $\mathfrak{g}$ is called square integrable (SQ.I.) if $i(\mathfrak{g}) = \dim Z(\mathfrak{g})$. (In the nilpotent case such Lie algebras are precisely the Lie algebras of simply connected nilpotent Lie groups having square integrable representations [MW, Theorem 3]). For instance any Heisenberg Lie algebra is square integrable.\\
Choose a basis $x_1, \ldots, x_t, x_{t+1},\ldots, x_n$ such that $x_{t+1},\ldots, x_n$ is a basis of $Z(\mathfrak{g})$.\\
Then, $t = \dim \mathfrak{g} - \dim Z(\mathfrak{g}) = \dim \mathfrak{g} - i(\mathfrak{g})$, which is the rank of the matrix $([x_i, x_j])_{1\leq i,j\leq t}$.  By the above, its Pfaffian coincides with $p_{\mathfrak{g}}$ (up to a nonzero scalar).  Hence, $\deg p_{\mathfrak{g}} = (\dim \mathfrak{g} - i(\mathfrak{g}))/2 \geq 1$ and so $\mathfrak{g}$ is singular. In particular, any Frobenius Lie algebra $\mathfrak{g}$ (i.e. $i(\mathfrak{g}) = 0$ [O2, E]) is singular.\\
\ \\
{\bf 2.4 Commutative polarizations of $\mathfrak{g}$}\\
Put $c(\mathfrak{g}) = (\dim \mathfrak{g} + i(\mathfrak{g}))/2$.  Now, suppose $\mathfrak{g}$ admits a commutative Lie subalgebra $\mathfrak{h}$ such that $\dim \mathfrak{h} = c(\mathfrak{g})$, i.e. $\mathfrak{h}$ is a commutative polarization (notation: CP) with respect to any $\xi \in \mathfrak{g}_{\reg}^\ast$ [D5, 1.12].\\
These CP's occur frequently in the nilpotent case.  If in addition $\mathfrak{h}$ is an ideal of $\mathfrak{g}$ then we call $\mathfrak{h}$ a CP-ideal (notation: CPI).  If a solvable Lie algebra $\mathfrak{g}$ admits a CP then it also admits a CPI [EO, Theorem 4.1]. \\
CP's play a special role in the construction of irreducible representations of $U(\mathfrak{g})$ and their kernels, the primitive ideals [EO, pp. 140-141]. However, for our purposes the following is more useful.\\
\ \\
{\bf Theorem 7.} [O3, Theorems 7,14; Corollary 16]\\
Let $\mathfrak{h}$ be a commutative Lie subalgebra of $\mathfrak{g}$.  Then,
\begin{itemize}
\item[(i)] $\dim \mathfrak{h} \leq c(\mathfrak{g})$
\item[(ii)] If $\mathfrak{h}$ is a CP of $\mathfrak{g}$ then $R(\mathfrak{h})$ is a maximal Poisson commutative subfield of $R(\mathfrak{g})$ (equivalently: $D(\mathfrak{h})$ is a maximal subfield of $D(\mathfrak{g})$).\\
The converse holds if $\mathfrak{g}$ is algebraic.\\
\ \\
\end{itemize}

{\bf 2.5 The Frobenius semi-radical $F(\mathfrak{g})$}\\
Put $F(\mathfrak{g}) = \sum\limits_{\xi \in \mathfrak{g}_{\reg}^{\ast}} \mathfrak{g}(\xi)$.  This is a characteristic ideal of $\mathfrak{g}$ containing $Z(\mathfrak{g})$ and for which $F(F(\mathfrak{g})) = F(\mathfrak{g})$.  It can also be characterized as follows: $R(\mathfrak{g})^{\mathfrak{g}} \subset R(F(\mathfrak{g}))$ and if $\mathfrak{g}$ is algebraic then $F(\mathfrak{g})$ is the smallest Lie subalgebra of $\mathfrak{g}$ with this property.  Similar results hold in $D(\mathfrak{g})$  [O4, Proposition 2.4, Theorem 2.5] (see also [Ce2]).\\
As a special case we have the following:\\
\ \\
{\bf Remark 8.} $Y(\mathfrak{g}) \subset S(F(\mathfrak{g}))$ (respectively $Z(U(\mathfrak{g})) \subset U(F(\mathfrak{g})))$ and $F(\mathfrak{g})$ is the smallest Lie subalgebra of $\mathfrak{g}$ with this property in case $\mathfrak{g}$ is an algebraic Lie algebra without proper semi-invariants.\\
\ \\
In case $\mathfrak{g}$ is square integrable we notice that $F(\mathfrak{g}) = Z(\mathfrak{g})$ (since $\mathfrak{g}(\xi) = Z(\mathfrak{g})$ for all regular $\xi \in \mathfrak{g}^\ast$) which forces $R(\mathfrak{g})^{\mathfrak{g}} = R(Z(\mathfrak{g}))$.  See also [O4, p. 283; DNOW, p. 323].  In particular, $Y(\mathfrak{g}) = S(Z(\mathfrak{g}))$, which is a polynomial algebra.\\
\ \\
If $\mathfrak{g}$ admits a CP $\mathfrak{h}$ then $F(\mathfrak{g})$ is commutative (since $F(\mathfrak{g}) \subset \mathfrak{h})$.  We will examine in Theorem 22 and Remark 23 when the converse holds. Clearly,
\begin{center}
$F(\mathfrak{g}) = 0$ if and only if $\mathfrak{g}$ is Frobenius\\
\end{center}
For this reason $F(\mathfrak{g})$ is called the Frobenius semi-radical of $\mathfrak{g}$.  At the other end of the spectrum we have the Lie algebras for which $F(\mathfrak{g}) = \mathfrak{g}$, which we call quasi-quadratic.  These are unimodular and they do not possess any proper semi-invariants.  They form a large class, which include all quadratic Lie algebras (and hence all abelian and semi-simple Lie algebras) [O4].\\
\ \\
{\bf 2.6 The sum rule and its converse}\\
{\bf Theorem 9.} [OV, Proposition 1.4]\\
Assume that $Y(\mathfrak{g})$ is freely generated by homogeneous elements $f_1,\ldots, f_r$.\\
If $\mathfrak{g}$ has no proper semi-invariants $(\bullet)$ then:
$$\sum\limits_{i=1}^r \deg f_i = c(\mathfrak{g}) - \deg p_{\mathfrak{g}}$$
where $p_{\mathfrak{g}}$ is the fundamental semi-invariant of $\mathfrak{g}$.\\
\ \\
{\bf Remark 10.} Recently Joseph and Shafrir [JS] were able to extend this sum rule by replacing condition $(\bullet)$ by: $\trdeg_k\ Y(\mathfrak{g}) = i(\mathfrak{g})$, $\mathfrak{g}$ is unimodular and $p_{\mathfrak{g}}$ is an invariant.\\ 
\ \\
{\bf Proposition 11.} Let $B$ be the Borel subalgebra of a simple Lie algebra $\mathfrak{g}$ of type $A_n, n\geq 2$, and let $N$ be its nilradical.  Then,
\begin{itemize}
\item[(i)] There are homogeneous, algebraically independent generators $f_1,\ldots, f_m$ of $Y(N)$ with $\deg f_i = i$ for all $i:1,\ldots, m$ and $i(N) = m = \left[\frac{1}{2}(n+1)\right]$.
\item[(ii)] $\deg p_N = \frac{1}{2} t(t+1)$ where $t = \left[\frac{1}{2} n\right]$.  In particular $N$ is singular.
\item[(iii)] If $n$ is odd then $F(N)$ is a CP of $N$.  If $n$ is even then $F(N)$ is the intersection of two CP's of $N$.
\item[(iv)] Suppose $n > 2$.  Then $c(N) - c(F(N)) < \deg p_N$.\\
\end{itemize}

{\bf Proof.} (i) This is due to Dixmier [D3] since $N$ can be identified with the Lie algebra of strictly lower triangular $(n + 1) \times (n + 1)$ matrices.  Clearly, $\dim\ N = \ds\frac{1}{2} n(n+1)$ and $\sum\limits_{i=1}^m \deg f_i = \sum\limits_{i=1}^m i = \ds\frac{1}{2} m(m+1)$
\begin{itemize}
\item[(ii)] There are two cases to distinguish:\\
1) $n = 2t$. \hs Then $i(N) = m = t$ and
$$c(N) = \frac{1}{2} (\dim N + i(N)) = \frac{1}{2}\left(t(2t + 1) + t\right) = t(t+1)$$
$$\deg p_N = c(N) - \sum\limits_{i=1}^m \deg f_i = t(t+1) - \frac{1}{2} t(t+1) = \frac{1}{2} t(t+1)$$
2) $n = 2t+1$.\ Then $i(N) = m = t + 1$.
$$c(N) = \frac{1}{2} (\dim\ N + i(N)) = \frac{1}{2} \left((t+1)(2t+1) + t+1\right) = (t+1)^2$$
$$\deg p_N = c(N) - \sum\limits_{i=1}^m \deg f_i = (t+1)^2 - \frac{1}{2} (t+1) (t+2) = \frac{1}{2} t(t+1)$$
\end{itemize}
Consequently, $N$ is singular by Remark 5.
\begin{itemize}
\item[(iii)] See [O4, Theorem 4.1].
\item[(iv)] This is obvious if $n$ is odd since then $c(N) = c(F(N))$.  If $n$ is even, say $n = 2t$, then $F(N)$ is commutative by (iii).  Hence, $c(F(N)) = \dim F(N) = t^2$ by [O4, Theorem 4.1].  It follows that:
$$c(N) - c(F(N)) = t(t+1) - t^2 = t < \frac{1}{2} t(t+1) = \deg p_N.$$
\end{itemize}

{\bf Remark 12.} In case $\mathfrak{g}$ is simple of type $C_n$ then $N$ is coregular [J2] and nonsingular since $\codim\ N_{\sing}^\ast = 2$ [J6, Lemma 2.6.17].  Furthermore, $N$ admits a CP [EO, Theorem 6.1] which is unique since it coincides with $F(N)$.\\
\ \\
Each of the following provides a test for coregularity.  The first two can be found in [OV, Corollary 1.2].  We add a minor extension to the first, while the third one will be shown in the more general setting of Theorem 26.\\
\ \\
{\bf Theorem 13.} Let $\mathfrak{g}$ be a nonabelian, coregular Lie algebra without proper semi-invariants.  Let $f_1,\ldots, f_r$ be homogeneous, algebraically independent generators of $Y(\mathfrak{g})$.  Then, 
\begin{itemize}
\item[(1)] $3i(\mathfrak{g}) + 2\deg p_{\mathfrak{g}} \leq \dim \mathfrak{g} + 2 \dim Z(\mathfrak{g})$\\
Moreover, equality occurs if and only if $\deg f_i \leq 2$, $i:1,\ldots, r$.
\item[(2)] $\codim \ \mathfrak{g}_{\sing}^{\ast} \leq 3$
\item[(3)] Suppose in addition that $\mathfrak{g}$ is algebraic.  If $\mathfrak{g}$ admits a CP then $\codim\ \mathfrak{g}_{\sing}^{\ast} \leq 2$.
\end{itemize}

The following is a converse of the sum rule.  It extends [PPY, Theorem 1.2], [P, Theorem 1.2], where the nonsingular algebraic case was treated.  It provides a powerful tool for proving the polynomiality of $Y(\mathfrak{g})$.\\
\ \\
{\bf Theorem 14.} [JS, 5.7]\\
Assume $\trdeg_k Y(\mathfrak{g}) = i(\mathfrak{g})$.  Let $f_1,\ldots, f_r \in Y(\mathfrak{g})$, $r = i(\mathfrak{g})$, be homogeneous, algebraically independent polynomials such that
$$\sum\limits_{i=1}^r \deg f_i \leq c(\mathfrak{g}) - \deg p_{\mathfrak{g}}$$
Then $Y(\mathfrak{g})$ is freely generated by $f_1,\ldots, f_r$.\\
\ \\
{\bf 2.7 Maximal Poisson commutative subalgebras of $S(\mathfrak{g})$}\\
Let $A$ be a Poisson commutative subalgebra of $S(\mathfrak{g})$.  Then it is well-known that $\trdeg_k(A) \leq c(\mathfrak{g})$.  $A$ is called complete if equality holds and strongly complete if it is also a maximal Poisson commutative subalgebra.  Denote by $Q(A)$ its quotient field and by $A'$ its Poisson commutant in $S(\mathfrak{g})$.\\
\ \\
{\bf Proposition 15.}  Let $\mathfrak{g}$ be algebraic and let $A$ be a complete, Poisson commutative subalgebra of $S(\mathfrak{g})$.  Then,
\begin{itemize}
\item[(i)] [JS, 7.1] $A'$ is Poisson commutative and strongly complete.
\item[(ii)] [PY, 2.1] $A$ is strongly complete if and only if $A$ is algebraically closed in $S(\mathfrak{g})$ (i.e. if $f \in S(\mathfrak{g})$ is algebraic over $A$ then $f \in A$).
\end{itemize}
The following provides a useful criterion in order to verify the condition of (ii).  It is the characteristic zero version by Panyushev, Premet and Yakimova [PPY, Theorem 1.1], [P, Theorem 1.5] of a result by Skryabin [Sk, Theorem 5.4] in positive characteristic.\\
\ \\
{\bf Theorem 16.} Let $A$ be a subalgebra of $S(\mathfrak{g})$ generated by homogeneous elements $f_1,\dots, f_r \in A$.  Consider the Jacobian locus:
$$J(f_1,\ldots, f_r) = \{\xi \in \mathfrak{g}^{\ast}\mid d_\xi f_1,\ldots, d_\xi f_r\ \mbox{are linearly dependent}\}$$
If $\codim \ J(f_1,\ldots, f_r) \geq 2$ in $\mathfrak{g}^\ast$, then $A$ is  algebraically closed in $S(\mathfrak{g})$.\\
\ \\
{\bf Remark 17.} Let $M$ be a maximal Poisson commutative subalgebra of $S(\mathfrak{g})$.  Then,
\begin{itemize}
\item[(i)] Suppose $f$, $g$ are nonzero elements of $S(\mathfrak{g})$.  If $fg $ and $g$ belong to $M$, then so does $f$.
\item[(ii)] $M$ is not necessarily complete and $Q(M)$ is not necessarily a maximal Poisson commutative subfield of $R(\mathfrak{g})$.
\end{itemize}

{\bf Proof.} (i) Clearly $Q(M) \cap S (\mathfrak{g})$ is a Poisson commutative subalgebra of $S(\mathfrak{g})$ containing $M$ and therefore coincides with $M$.  Hence $f = (fg) g^{-1} \in Q(M) \cap S(\mathfrak{g}) = M$.\\
(ii) Consider the 3-dimensional, algebraic Lie algebra $\mathfrak{g}$ with basis $x, y, z$, with nonzero brackets $[x,y] = y$ and $[x ,z] = z$. Clearly, $i(\mathfrak{g}) = 1$, and $c(\mathfrak{g}) = 2$.  One easily verifies that $M = k[x]$ is a maximal Poisson commutative subalgebra of $S(\mathfrak{g})$.  However, it is not complete as $\trdeg_k M = 1 < 2 = c(\mathfrak{g})$.  Also, $Q(M) = k(x)$ is not a maximal Poisson commutative subfield of $R(\mathfrak{g})$ since it does not contain $y/z$ which Poisson commutes with $x$.\\
\ \\
{\bf Proposition 18.} Let $A$ be a subalgebra of $S(\mathfrak{g})$ such that $Q(A)$ is a maximal Poisson commutative subfield of $R(\mathfrak{g})$.  Then,
\begin{itemize}
\item[(i)] $A'=Q(A) \cap S(\mathfrak{g})$, which is a maximal Poisson commutative subalgebra of $S(\mathfrak{g})$ and $Q(A) = Q(A')$.  Similar results hold in $U(\mathfrak{g})$.
\item[(ii)] If $\{\mathfrak{g},A\} \subset A$ then $\{\mathfrak{g},A'\}\subset A'$.
\item[(iii)] Suppose $B$ is a commutative subalgebra of $U(\mathfrak{g})$ such that the symmetrization $s: A \rightarrow B$ is an associative algebra isomorphism, then the same is true of
$$s:A_c \cap S(\mathfrak{g}) \rightarrow B_c \cap U(\mathfrak{g})$$
for any nonzero $c \in A \cap Z(\mathfrak{g})$, where $A_c$ (resp. $B_c$) is the localization of $A$ (resp. $B$) at $c$.
\item[(iv)] Suppose that $A$ is a graded subalgebra of $S(\mathfrak{g})$ and that $Q(A) \cap S(\mathfrak{g}) = A_c \cap S(\mathfrak{g})$ for a suitable $c \in A \cap Z(\mathfrak{g})$. Then $Q(B) \cap U(\mathfrak{g}) = B_c \cap U(\mathfrak{g})$.
\end{itemize}

{\bf Proof.}
\begin{itemize}
\item[(i)] Put $M = Q(A) \cap S(\mathfrak{g})$.  Clearly $M \subset A'$ since $M$ Poisson commutes with $A$.  On the other hand, take $x \in A'$, i.e. $\{x,A\} = 0$, which implies that $\{x, Q(A)\} = 0$.  By the maximality of $Q(A)$ we obtain that $x \in Q(A) \cap S(\mathfrak{g}) = M$.  Consequently, $M = A'$.  Next, we have to show that $M' = M$.  Indeed, $M \subset M'$ by the Poisson commutativity of $M$, while $A \subset M$ implies that $M' \subset A' = M$.  Finally, $Q(A) = Q(A')$ follows easily from the fact that $A \subset A' \subset Q(A)$.
\item[(ii)] This is straightforward and well-known.
\item[(iii)] This can be easily deduced from the fact that $s(cx) = cs(x)$ for all $x \in A$ and nonzero $c \in A \cap Z(\mathfrak{g})$.
\item[(iv)] Let $(U_q)_{q \geq 0}$ be the natural increasing filtration of $U(\mathfrak{g})$.  The associated graded algebra $gr(U(\mathfrak{g}))$ can be identified with $S(\mathfrak{g})$. 
\end{itemize}
The elements $u \in U_q\backslash U_{q-1}$ are said to be of degree $q$ and $[u] = u\ mod\ U_{q-1}$ is called the leading term of $u$. Obviously, $B_c \cap U(\mathfrak{g}) \subset Q(B) \cap U(\mathfrak{g})$.  Next, we take a nonzero $w \in Q(B) \cap U(\mathfrak{g})$. By induction on $n = deg\ w$ we show that $w \in B_c \cap U(\mathfrak{g})$.   Since $B = s(A)$ we can find $x, y \in A$ such that $w = s(x) s(y)^{-1}$.  Let $x = x_p + \ldots + x_0$, $x_p \neq 0$, and $y = y_q + \ldots + y_0$, $y_q \neq 0$, be the decomposition of $x$ and $y$ into homogeneous components, which belong to $A$ by assumption.  In particular, $[s(x)] = x_p \in A$ and $[s(y)] = y_q \in A$.  Now we observe that
$$x_p = [s(x)] = [ws(y)] = [w] [s(y)] = [w] y_q$$
Hence, $ [w] = x_p y_q^{-1} \in Q(A) \cap S(\mathfrak{g}) = A_c \cap S(\mathfrak{g})$ and so\\
$s([w]) \in B_c\cap U(\mathfrak{g}) \subset Q(B) \cap U(\mathfrak{g})$.\\
Clearly, $w-s([w]) \in Q(B)\ \cap\ U(\mathfrak{g})$ and its degree is strictly less than $n$.  By the induction hypothesis, $w-s([w]) \in B_c\ \cap\ U(\mathfrak{g})$.  Finally, $w = (w - s([w])) + s([w]) \in B_c\ \cap\ U(\mathfrak{g})$.\\
\ \\
{\bf Examples 19.}
\begin{itemize}
\item[(1)] Suppose $\mathfrak{g}$ admits a CP $\mathfrak{h} \subset \mathfrak{g}$.  Then $S(\mathfrak{h})$ is a polynomial, strongly complete subalgebra of $S(\mathfrak{g})$ (Indeed, by Theorem 7 $R(\mathfrak{h})$ is a maximal Poisson commutative subfield of $R(\mathfrak{g})$, $S(\mathfrak{h}) = R(\mathfrak{h}) \cap S(\mathfrak{g})$ and $\trdeg_k S(\mathfrak{h}) = \dim \mathfrak{h} = c(\mathfrak{g})$).  If $\mathfrak{h}$ is a CPI of $\mathfrak{g}$ then $\{\mathfrak{g},S(\mathfrak{h})\} \subset S(\mathfrak{h})$.  Similar results hold for $U(\mathfrak{h})$.
\item[(2)] Let $\mathfrak{g}$ be an $n$-dimensional, nilpotent Lie algebra. Choose an increasing sequence of ideals $\mathfrak{g}_i$ of $\mathfrak{g}$ such that $\dim \mathfrak{g}_i = i$, $i: 0,\ldots, n$.  Following Mich\`ele Vergne [V] we let $A$ (respectively $B$) be the subalgebra of $S(\mathfrak{g})$ (resp. $U(\mathfrak{g})$) generated by the union of $Y(\mathfrak{g}_i)$ (resp. $Z(U(\mathfrak{g}_i)))$.  Then $\{\mathfrak{g}, A\} \subset A$, $[\mathfrak{g},B] \subset B$ and the symmetrization $s : A \rightarrow B$ is an associative algebra isomorphism.\\
The quotient field $Q(A)$ (resp. $Q(B)$) is a maximal subfield of the Poisson field $R(\mathfrak{g})$ (resp. of $D(\mathfrak{g})$), which is a purely transcendental extension of $k$ of degree $c(\mathfrak{g})$.  Since $\mathfrak{g}$ is nilpotent $U(\mathfrak{g})$ satisfies the Gelfand-Kirillov conjecture [GK].  In [V] the analogue of this is shown for the Poisson algebra $S(\mathfrak{g})$.\
By Proposition 18, $M = Q(A) \cap S(\mathfrak{g}) = A'$ is strongly complete, $\{\mathfrak{g},M\} \subset M$ and similarly $N = Q(B) \cap U(\mathfrak{g})$ is a maximal commutative subalgebra of $U(\mathfrak{g})$.\\
Moreover, if $\mathfrak{g}$ is indecomposable of dimension at most seven we will see in section 3 that $M$ is polynomial and that $M = A_c\cap S(\mathfrak{g})$ for some nonzero $c \in A \cap Z(\mathfrak{g})$.  Consequently, the symmetrization $s : M \rightarrow N$ is an associative algebra isomorphism by (iii) and (iv) of Proposition 18 as $A$ is a graded subalgebra of $S(\mathfrak{g})$.  In the sequel we will denote $A$ by $V(\mathfrak{g})$ for a fixed sequence of ideals.
\item[(3)] The first part of Proposition 18 can also be employed if $\mathfrak{g}$ is solvable, because then there exists a polynomial subalgebra $B \subset U(\mathfrak{g})$ such that $Q(B)$ is a maximal subfield of $D(\mathfrak{g})$ [N].  This result played an essential role in Joseph's proof of the Gelfand-Kirillov conjecture in the solvable case [J1].\\
\end{itemize}

{\bf 2.8 The Mishchenko-Fomenko algebras and the Bolsinov criterion}\\
Let $x_1,\ldots, x_n$ be a basis of $\mathfrak{g}$ and let $f \in S(\mathfrak{g})$ be a polynomial of degree $d$ in the $x_i$'s.  Let $\xi \in \mathfrak{g}^\ast$, $t \in k$ and put $a_i = \xi(x_i) \in k$, $i = 1,\ldots,n$.  Consider the expansion
$$f(x_1 + a_1t , \ldots, x_n + a_nt) = \sum\limits_{j=0}^d f_{\xi}^j (x_1,\ldots, x_n)t^j$$
which gives rise to $d$ polynomials $f_{\xi}^j$ in $S(\mathfrak{g})$, $j: 0,1,\ldots, d-1$, called the $\xi$-shifts of $f$, with $\deg f_\xi^j = d-j$.\\
Note that $f_{\xi}^d \in k$ and is therefore left out.  Also, $f_{\xi}^0 = f$ and $f_{\xi}^{d-1} = d_{\xi}f$.\\
Furthermore, if $f$ is homogeneous, then so are the $\xi$-shifts.\\
Denote by $Y_\xi(\mathfrak{g})$ the subalgebra of $S(\mathfrak{g})$ generated by the $\xi$-shifts of all $f \in Y(\mathfrak{g})$ (or equivalently of the generators of $Y(\mathfrak{g})$).  The subalgebras of the form $Y_\xi(\mathfrak{g})$ are called the Mishchenko-Fomenko algebras.  It was observed by Mishchenko and Fomenko that these subalgebras are Poisson commutative [MF] and hence $\trdeg_k Y_\xi (\mathfrak{g}) \leq c(\mathfrak{g})$.  They also showed that for $\mathfrak{g}$ semi-simple, $Y_\xi(\mathfrak{g})$ is complete (i.e. that equality occurs) for all $\xi \in \mathfrak{g}_{\reg}^\ast$ [MF].\\
A general criterion by Bolsinov asserts that for any $\mathfrak{g}$ and $\xi \in \mathfrak{g}_{\reg}^\ast$ that $Y_\xi(\mathfrak{g})$ is complete if and only if $\mathfrak{g}$ is nonsingular [Bo1, Bo2].\\
However, one must here include the condition that $\trdeg_kY(\mathfrak{g}) = i(\mathfrak{g})$ as pointed out by Panyushev and Yakimova [PY, Theorem 2.3]. We now provide the following.\\
\ \\
{\bf Counterexample 20} (to Bolsinov's assertion)\\
It is example (5.8) of [DDV, p. 322], namely the 8-dimensional solvable Lie algebra $\mathfrak{g}$ with basis $x_0, x_1,\ldots, x_7$ over $k$ with nonzero brackets:\\
$[x_0, x_1] = 5x_1$, $[x_0, x_2] = 10x_2$, $[x_0, x_3] = -13x_3$, $[x_0, x_4] = -8x_4$,\\
$[x_0, x_5] = -3x_5$, $[x_0, x_6] = 2x_6$, $[x_0, x_7] = 7x_7$, $[x_1, x_3] = x_4$,\\
$[x_1, x_4] = x_5$, $[x_1, x_5] = x_6$, $[x_1,x_6] = x_7$, $[x_2, x_3] = x_5$, $[x_2, x_4] = x_6$, $[x_2, x_5] = x_7$.\\
Then, $\mathfrak{g}$ is algebraic and unimodular with $i(\mathfrak{g}) = 2$.  Also, $Y(\mathfrak{g}) = k$ (and so $\trdeg_k Y(\mathfrak{g}) = 0 < i(\mathfrak{g}))$ which implies that $Y_{\xi} (\mathfrak{g}) = k$ for any $\xi \in \mathfrak{g}_{\reg}^\ast$.  Therefore, $\trdeg_k Y_{\xi}(\mathfrak{g}) = 0 < 5 = c(\mathfrak{g})$, while $\mathfrak{g}$ is nonsingular (in fact $\codim\ \mathfrak{g}_{\sing}^\ast = 3$).  Also note that the extended sum rule (Remark 10) is not valid in this situation.\\
Finally, $\mathfrak{h} = < x_3, x_4, x_5, x_6, x_7>$ is a CPI for $\mathfrak{g}$\\
\ \\
Recently, Joseph and Shafrir constructed an open subset $\mathfrak{g}_{w\reg}^\ast \subset \mathfrak{g}^\ast$ (for details see [JS, 7.2]) for which $\mathfrak{g}_{\reg}^\ast \subset \mathfrak{g}_{w\reg}^\ast$ with equality if and only if $\mathfrak{g}$ is nonsingular.\\
Furthermore, they obtained the following extension (and correction) of the Bolsinov criterion.\\
\ \\
{\bf Theorem 21.} Suppose $\trdeg_k Y(\mathfrak{g}) = i(\mathfrak{g})$.  Then, 
$$\trdeg_k Y_\xi(\mathfrak{g}) \leq c(\mathfrak{g}) - \deg p_\mathfrak{g}$$
with equality if and only if $\xi \in \mathfrak{g}_{w\reg}^\ast$.\\
\ \\
Now, we want to compare $c(F(\mathfrak{g}))$ with $c(\mathfrak{g})$.  For instance they coincide if $F(\mathfrak{g})$ is of codimension one in $\mathfrak{g}$ by [EO, (4) of Proposition 1.6].  At the same time we will examine when the commutativity of $F(\mathfrak{g})$ guarantees the existence of a CP in $\mathfrak{g}$ (the converse is always valid).\\
\ \\
{\bf Theorem 22.}
\begin{itemize}
\item[(1)] $c(\mathfrak{h}) \leq c(\mathfrak{g})$ for any Lie subalgebra $\mathfrak{h}$ of $\mathfrak{g}$.
\item[(2)] Assume that $\trdeg_k Y(\mathfrak{g}) = i(\mathfrak{g})$.  Then,
$$c(\mathfrak{g}) - c(F(\mathfrak{g})) \leq \deg p_g\insp (\ast\ast)$$
So, $c(\mathfrak{g}) = c(F(\mathfrak{g}))$ if $\mathfrak{g}$ is nonsingular (by Remark 5).  In particular, any complete Poisson commutative subalgebra of $S(F(\mathfrak{g}))$ is also complete in $S(\mathfrak{g})$.  If in addition $F(\mathfrak{g})$ is commutative then $F(\mathfrak{g})$ is the only CP of $\mathfrak{g}$.
\item[(3)] Strict inequality can occur in $(\ast\ast)$.
\end{itemize}

{\bf Proof.} 
\begin{itemize}
\item[(1)] Fix $\xi \in g_{\reg}^\ast$.  Then, $B_\xi(x,y) = \xi([x,y])$ is an alternating bilinear form on $\mathfrak{g}$.  Consider $\xi\mid_{\mathfrak{h}} \in \mathfrak{h}^\ast$. By [D5, 1.12.17] there exists a solvable Lie subalgebra $P$ of $\mathfrak{h}$ of dimension $\frac{1}{2}(\dim \mathfrak{h} + i(\mathfrak{h})) = c(\mathfrak{h})$ which is subordinate to $\xi\mid_{\mathfrak{h}}$ and hence also to $\xi$.  Because of the latter we may conclude by [D5, 1.12.1] that
$$c(\mathfrak{h}) = \dim P \leq \frac{1}{2} (\dim \mathfrak{g} + \dim \mathfrak{g}(\xi)) = \frac{1}{2} (\dim \mathfrak{g} + i(\mathfrak{g})) = c(\mathfrak{g})$$
\item[(2)] Let $x_1, \ldots, x_s$, $x_{s+1},\ldots, x_n$ be a basis of $\mathfrak{g}$ such that $x_1,\ldots, x_s$ is a basis of $F(\mathfrak{g})$.  Fix $\xi \in \mathfrak{g}_{\reg}^\ast$.  $Y_\xi(\mathfrak{g})$ is generated by the $\xi$-shifts $f_{\xi}^j$ of all $f \in Y(\mathfrak{g})$. By Remark 8, $Y(\mathfrak{g}) \subset S(F(\mathfrak{g}))$.  So, any§ $f \in Y(\mathfrak{g})$ depends only on $x_1,\ldots, x_s$ and so does each of its $\xi$-shifts.  Therefore $Y_\xi(\mathfrak{g}) \subset S(F(\mathfrak{g}))$.  From the Poisson commutativity of $Y_\xi(\mathfrak{g})$ and by using Theorem 21 we deduce:
$$c(\mathfrak{g}) - \deg p_{\mathfrak{g}} = \trdeg_k Y_\xi(\mathfrak{g}) \leq c(F(\mathfrak{g}))$$
which implies $(\ast\ast)$.\\
Suppose that $\mathfrak{g}$ is nonsingular and that $F(\mathfrak{g})$ is commutative.  Then $\dim F(\mathfrak{g}) = c(F(\mathfrak{g})) = c(\mathfrak{g})$, i.e. $F(\mathfrak{g})$ is a CP of $\mathfrak{g}$.  On the other hand, let $\mathfrak{h}$ be a CP of $\mathfrak{g}$.  Then we know that $F(\mathfrak{g}) \subset \mathfrak{h}$.  Having the same dimension we conclude that $F(\mathfrak{g}) = \mathfrak{h}$.
\item[(3)] Indeed, see for example (iv) of Proposition 11.
\end{itemize}

{\bf Remark 23.} (See also [EO, Proposition 2.2])\\
For any nilpotent Lie algebra $\mathfrak{g}$ of dimension at most 7 we have that:
\begin{center}
$\mathfrak{g}$ admits a CP \ \ $\Leftrightarrow \ \ F(\mathfrak{g})$ is commutative
\end{center}
In higher dimension this is no loànger valid.  Indeed in [EO, 3.1] examples of square integrable, nilpotent Lie algebras $\mathfrak{g}$ without CP's are given.  In particular, $F(\mathfrak{g}) = Z(§\mathfrak{g})$ is commutative.\\
\ \\
{\bf Proof.}  In view of the previous theorem it suffices to verify that $F(\mathfrak{g})$ is not commutative for any singular indecomposable nilpotent Lie algebra $\mathfrak{g}$, $\dim \mathfrak{g} \leq 7$, without CP's, i.e. for the Lie algebras 21, 69, 77, 100, 101 (see sections 3.2 and 3.3).\\
\ \\
The following is an abbreviated version of an important extension by Panyushev and Yakimova [PY] of a result by Tarasov [T] in the semi-simple case.\\
\ \\
{\bf Theorem 24.} Suppose that $\mathfrak{g}$ is algebraic for which:
\begin{itemize}
\item[(i)] $Y(\mathfrak{g})$ is freely generated by homogeneous polynomials $f_1,\ldots, f_r$, where $r = i(\mathfrak{g})$, such that $\sum\limits_{i=1}^r \deg f_i = c(\mathfrak{g})$.
\item[(ii)] $\codim\ \mathfrak{g}_{\sing}^\ast \geq 3$
\end{itemize}
Then, for any $\xi \in \mathfrak{g}_{\reg}^\ast$, $Y_\xi(\mathfrak{g})$ is a polynomial, strongly complete subalgebra of $S(\mathfrak{g})$.\\
\ \\
{\bf Example 25.} Let $\mathfrak{g}$ be the diamond Lie algebra over $k$, i.e. the 4-dimensional solvable Lie algebra with basis $t$, $x$, $y$, $z$ with nonvanishing brackets\\
$[t,x] = -x$, $[t,y] = y$, $[x,y] =z$.\\
$\mathfrak{g}$ is algebraic and it is easy to verify that $i(\mathfrak{g}) = 2$ and $\codim \ \mathfrak{g}_{\sing}^\ast = 3$.\\
$Y(\mathfrak{g})$ is freely generated by the homogeneous polynomials $xy - tz$ and $z$ of which the sum of the degrees equals $3 = c(\mathfrak{g})$.\\
Fix $\xi = z^\ast \in \mathfrak{g}_{\reg}^\ast$.  Then,
$$Y_{\xi} (\mathfrak{g}) = k[xy - tz, t, z] = k[xy, t, z]$$
which is a polynomial, strongly complete subalgebra of $S(\mathfrak{g})$ by the previous theorem.  Note that $xy \in Y_\xi(\mathfrak{g})$ but $x \notin Y_\xi(\mathfrak{g})$ and $y \notin Y_\xi(\mathfrak{g})$ (compare this with (i) of Remark 17).  Also, $\{\mathfrak{g}, Y_\xi(\mathfrak{g})\}\not\subset Y_\xi(\mathfrak{g})$ (since $\{x,xy\} = xz \notin Y_\xi(\mathfrak{g}))$.\\
\ \\
On the other hand, take $\eta = x^\ast \in \mathfrak{g}_{\reg}^\ast$.  Then, $Y_\eta (\mathfrak{g}) = k[xy - tz, y, z]$, which is also a polynomial, strongly complete subalgebra of $S(\mathfrak{g})$ and this time $\{\mathfrak{g},Y_\eta(\mathfrak{g})\} \subset Y_\eta(\mathfrak{g})$.\\
\ \\
{\bf Theorem 26.}  Let $\mathfrak{g}$ be a nonabelian unimodular, coregular Lie algebra such that $\trdeg_k Y(\mathfrak{g}) = i(\mathfrak{g})$.  If $\mathfrak{g}$ admits a CP $\mathfrak{h}$ then, for any $\xi \in \mathfrak{g}_{\reg}^\ast$, $Y_\xi(\mathfrak{g})$ is not a maximal Poisson commutative subalgebra of $S(\mathfrak{g})$.  If in addition $\mathfrak{g}$ is algebraic, then $\codim\ \mathfrak{g}_{\sing}^\ast \leq 2$.\\
\ \\
{\bf Remark.} This theorem is no longer valid if the condition $trdeg_k Y(\mathfrak{g}) = i(\mathfrak{g})$ is removed.  See counterexample 20.\\
\ \\
{\bf Proof.}  Suppose that $Y_\xi(\mathfrak{g})$ is a maximal Poisson commutative subalgebra for some $\xi \in \mathfrak{g}_{\reg}^\ast$.  Let $h_1,\ldots, h_m$, $m = c(\mathfrak{g})$, be a basis of $\mathfrak{h}$.  We know that $F(\mathfrak{g}) \subset \mathfrak{h}$.  As seen in the proof of Theorem 22,
$$Y_\xi(\mathfrak{g}) \subset S(F(\mathfrak{g})) \subset S(\mathfrak{h})$$
where $S(\mathfrak{h})$ is a Poisson commutative subalgebra of $S(\mathfrak{g})$.  It follows that $Y_\xi(\mathfrak{g}) = S(\mathfrak{h})$ by the maximality of $Y_\xi(\mathfrak{g})$.  By Theorem 21
$$c(\mathfrak{g}) - \deg p_{\mathfrak{g}} = \trdeg_k Y_\xi(\mathfrak{g}) = \trdeg_k S(\mathfrak{h}) = c(\mathfrak{g})$$
Hence $\deg p_{\mathfrak{g}} = 0$.\\
By assumption $Y(\mathfrak{g})$ is freely generated by homogeneous polynomials $f_1,\ldots, f_r$, where $r = i(\mathfrak{g})$.  Then the extended sum rule (Remark 10) asserts that
$$\sum\limits_{i=1}^r \deg f_i = c(\mathfrak{g})$$
$Y_\xi(\mathfrak{g})$ is generated by the $\xi$-shifts $f_i^j$, $i:1,\ldots, r$; $j:0,\ldots, \deg f_i-1$ of the $f_i$'s. Clearly their number is $\sum\limits_{i=1}^r \deg f_i = c(\mathfrak{g})$, which is precisely $\trdeg_k Y_\xi (\mathfrak{g})$.  Therefore $Y_\xi(\mathfrak{g})$ is freely generated by these $\xi$-shifts.  From the equality
$$k[f_i ^j \mid i:1,\ldots, r; \ j:0,\ldots, , \deg f_i-1] = k[h_1,\ldots, h_m]$$
we obtain that these two polynomial rings have the same Gorenstein invariant [OV, Example 5.6], i.e.
$$\sum\limits_{i,j} \deg f_i ^j = \sum\limits_{i=1}^m \deg h_i = m = c(\mathfrak{g})$$
Since the sum on the left hand side consists of $c(\mathfrak{g})$ positive integers, we deduce that $\deg\ f_i^j = 1$ for all $i$, $j$.  In particular, $\deg f_i = 1$, $i:1,\ldots, r$.  This implies that $Y(\mathfrak{g}) = Y_\xi(\mathfrak{g})$, hence $i(\mathfrak{g}) = c(\mathfrak{g}) = \frac{1}{2} (\dim \mathfrak{g} + i(\mathfrak{g}))$, i.e. $i(\mathfrak{g}) = \dim \mathfrak{g}$. Consequently, $\mathfrak{g}$ is abelian, contradicting our assumption.  Finally, suppose in addition that $\mathfrak{g}$ is algebraic, then we conclude that $\codim \ \mathfrak{g}_{\sing}^\ast \leq 2$ by theorem 24.\\
\ \\
{\bf Example 27: The standard filiform Lie algebra $\mathfrak{g}(n)$}\\
$\mathfrak{g}(n)$ is the nilpotent Lie algebra with basis $x_1,\ldots, x_n$ and nonzero brackets:
$$[x_1, x_i] = x_{i+1}, i:2,\ldots, n-1$$
$Y(\mathfrak{g}(n))$ is finitely generated [Ce4], but it is not always a polynomial algebra.  Dixmier verified by direct computation that $\mathfrak{g}(3)$ and $\mathfrak{g}(4)$ are coregular, but that $\mathfrak{g}(5)$ is not [D2].  In fact $\mathfrak{g} = \mathfrak{g}(n)$ is not coregular if $n \geq 5$ by (1) of Theorem 13 [OV, Example 1.7].  We now want to arrive at the same conclusion by invoking the previous theorem.  So, let us assume that $\mathfrak{g}(n)$ is coregular.  Clearly, $i(\mathfrak{g}(n)) = n-2$ and therefore $c(\mathfrak{g}(n)) = \frac{1}{2} (n+n-2) = n-1$.  Hence $\mathfrak{h} = <x_2,\ldots, x_n>$ is a CP of $\mathfrak{g}(n)$. One verifies that $\codim \ \mathfrak{g}_{\sing}^\ast = n-2$. Since $\mathfrak{g}(n)$ is nilpotent we may apply the previous theorem in order to conclude that $\codim \ \mathfrak{g}_{\sing}^\ast \leq 2$, i.e. $n \leq 4$.\\
\ \\
For the explicit generators of $Y(\mathfrak{g}(n))$, $n \leq 6$ we refer to [O5, section 5].  We now want to comment on $Y(\mathfrak{g}(7))$.  From the classical correspondence with the $SL(2)$-covariants of the binary quintic we know that $Y(\mathfrak{g}(7))$ is minimally generated by 23 elements, satisfying 168 relations [GY, pp. 131, 144].  An attempt to find these generators by using SINGULAR [GPS] failed because our computer ran out of memory.  However, later on we found them explicitly in an unpublished, carefully handwritten manuscript by Andr\'e Cerezo [Ce4], see also [Ce3].  The size of most of these generators is gigantic and their description takes about 20 pages.  For this reason we will only include in 3.2 (\# 159) the 5 algebraically independent homogeneous generators of the quotient field of $Y(\mathfrak{g}(7))$.\\
\ \\
{\bf Example 28.} (See also [JS, 8.3])\\
Let $\mathfrak{g}$ be the solvable Lie algebra with basis $x_0, x_1,\ldots, x_n$ and nonzero brackets $[x_0,x_i] = r_ix_i$, $i:1,\ldots,n$, where the $r_i$'s are (nonzero) rational numbers such that $\sum\limits_{i=1}^n r_i = 0$. Then,
\begin{center}
$\mathfrak{g}$ is coregular \ \ $\Leftrightarrow \ \ n = 2$
\end{center}
In that situation $Y(\mathfrak{g}) = k[x_1x_2]$.\\
\ \\
{\bf Proof.} Clearly $\mathfrak{g}$ is algebraic and unimodular with $i(\mathfrak{g}) = n-1$ and $n \geq 2$.  Hence $c(\mathfrak{g}) = \frac{1}{2} (n+1+n-1) = n$ and so $\mathfrak{h} = <x_1,\ldots,x_n>$ is a CP of $\mathfrak{g}$.  Also $\codim \ \mathfrak{g}_{\sing}^\ast = n$.  Next, we claim that $R(\mathfrak{g})^{\mathfrak{g}}$ is the quotient field $Q(Y(\mathfrak{g}))$ of $Y(\mathfrak{g})$.  Using the fact that each element of $R(\mathfrak{g})^{\mathfrak{g}}$ is the quotient of two semi-invariants of $S(\mathfrak{g})$ of the same weight, it is not difficult to verify that $R(\mathfrak{g})^{\mathfrak{g}}$ is generated as a field by elements of the form $z_m = x_1^{m_1} \ldots x_n^{m_n}$ where $m = (m_1,\ldots, m_n) \in \IZ^n$ such that $\sum\limits_{i=1}^n r_i m_i = 0$. In particular, $z_q = x_1 \ldots x_n \in Y(\mathfrak{g})$ where $q = (1,\ldots, 1)$.  Now consider $z_m$.  Choose a positive integer $t$ such that $m_i + t \geq 0$ for all $i : 1,\ldots, n$.  Then,
$$z_m(z_q)^t = z_{m+tq} = x_1^{m_1+t} \ldots x_n^{m_n+t} \in Y(\mathfrak{g})$$
hence, $z_m = z_{m+tq} (z_q)^{-t} \in Q(Y(\mathfrak{g}))$.\\
This establishes the claim.  Consequently,
$$\trdeg_k Y(\mathfrak{g}) = \trdeg_k R(\mathfrak{g})^{\mathfrak{g}} = i(\mathfrak{g})$$
by (1) of Theorem 1 since $\mathfrak{g}$ is algebraic.\\
Now, suppose $\mathfrak{g}$ is coregular.  Invoking the previous theorem, we deduce that\\
$n = \codim \ \mathfrak{g}_{\sing}^\ast \leq 2$ and thus $n = 2$.\\
The converse is obvious.\\
\ \\
{\bf Remark 29.}  One would expect Theorem 24 to be applicable for quadratic Lie algebras.  This is indeed the case for semi-simple Lie algebras [T], the diamond Lie algebra (Example 25) and the quadratic nilpotent Lie algebras of dimension at most 6 ($\mathfrak{g}_{5,4} (\# 7)$ and $g_{6,3} (\# 22)$ [O5,5]).  However this is not true for $\mathfrak{g}_{7,2.4} (\# 101)$ (See Example 32), which is the only indecomposable, quadratic nilpotent Lie algebra of dimension 7.  In fact, $\mathfrak{g}_{7,1.1(iii)} (\# 136)$ (See Example 33) is the only 7-dimensional indecomposable nilpotent Lie algebra which satisfies the conditions of Theorem 24.\\
\ \\
{\bf 2.9 A conjecture by Milovanov}\\
The following assertion raised by Milovanov, has been verified for a number of low dimensional Lie algebras by Korotkevich [K].\\
\ \\
{\bf Conjecture 30.} For any solvable Lie algebra $\mathfrak{g}$ there exists a complete Poisson commutative subalgebra $M$ of $S(\mathfrak{g})$ generated by elements of degree at most two.\\
\ \\
We thank Oksana Yakimova for informing us about this and for giving some helpful suggestions.  In the next section we will show this conjecture to be valid for all indecomposable (and hence also for all) at most 7-dimensional nilpotent Lie algebras.  This is obvious if $\mathfrak{g}$ admits a CP $\mathfrak{h}$, namely just take $M = S(\mathfrak{h})$, see (1) of Examples 19. Also, the subalgebra $M$ as constructed in (2) of Examples 19 has the required properties, except for the Lie algebras 134, 135, 136.  For these an alternative subalgebra is given, denoted by $M_1$, (see e.g. Example 33), satisfying the conjecture.\\
\ \\
{\bf \large 3. Indecomposable nilpotent Lie algebras of dimension $\leq 7$ ($k = \IC$)}\\
We will use the classification of Magnin [Ma1, Ma2, Ma3] and of Carles [Ca].  In view of the work already done in [O5,5], it suffices to consider all 7-dimensional indecomposable nilpotent Lie algebras $\mathfrak{g}$ for which $i(\mathfrak{g}) > \rank\ \mathfrak{g}$ (77 cases in total, of which $\mathfrak{g}_{7,0.4(\lambda)}$ (\# 84) is the only infinite family).\\
\ \\
{\bf 3.1 Procedure and examples}\\
Most of the calculations are done with MAPLE and with SINGULAR (if $\mathfrak{g}$ has an abelian ideal of codimension one).  We will primarily follow the method Dixmier used for dimension 5 [D2, pp. 322-330].  Therefore we recall the following special case of one of his results [D1, p. 333].\\
\ \\
{\bf Theorem 31.} Let $\mathfrak{g}$ be a nilpotent Lie algebra and let
$$0 = \mathfrak{g}_0 \subset \mathfrak{g}_1 \ldots \subset \mathfrak{g}_n = \mathfrak{g}$$
be a sequence of ideals of $\mathfrak{g}$ such that for each $j : 1,\ldots,n$, $\dim \mathfrak{g}_j = j$ and $[\mathfrak{g},\mathfrak{g}_j] \subset \mathfrak{g}_{j-i}$.  Choose $x_j \in \mathfrak{g}_j \backslash \mathfrak{g}_{j-1}$.  Suppose $j_1 < j_2 < \ldots < j_r$ are the indices $j \geq 1$ such that
$$S(\mathfrak{g}_{j-1}) \cap Y(\mathfrak{g}) \stackrel{\subset}{\neq} S(\mathfrak{g}_j) \cap Y(\mathfrak{g})$$
\begin{itemize}
\item[(1)] Then for each such $j$ there is a nonzero element $b_j \in S(\mathfrak{g}_{j-1}) \cap Y(\mathfrak{g})$ and $c_j \in S(\mathfrak{g}_{j-1})$ such that $a_j = b_j x_j + c_j \in S(\mathfrak{g}_j) \cap Y(\mathfrak{g})$.
\end{itemize}
In (2), (3), (4) $a_j$, $b_j$, $c_j$ are chosen to satisfy  (1).
\begin{itemize}
\item[(2)] $Y(\mathfrak{g}) \subset k[a_{j_1},\ldots, a_{j_r}, b_{j_1}^{-1},\ldots, b_{j_r}^{-1}]$
\item[(3)] $R(\mathfrak{g})^{\mathfrak{g}}$ is the quotient field of $Y(\mathfrak{g})$.  It is the field generated by $a_{j_1},\ldots, a_{j_r}$, which are algebraically independent over $k$.  In particular, $r = i(\mathfrak{g})$.
\item[(4)] $Y(\mathfrak{g}) \subset k[a_{j_1},\ldots, a_{j_r}, a^{-1}]$ for some nonzero $a \in k[a_{j_1},\ldots, a_{j_r}]$ (in our limited setting a can be taken in $Z(\mathfrak{g})$)
\end{itemize}

Using this we construct as a first step homogeneous, algebraically independent elements $f_1, \ldots, f_r \in Y(\mathfrak{g})$, which generate the quotient field $R(\mathfrak{g})^{\mathfrak{g}}$, $r = i(\mathfrak{g})$.  At this stage we recover the results of Romdhani [R], although quite a few cases are missing [G, pp. 146-149] and some minor corrections have to be made.  Next, we distinguish two major classes:\\
{\bf I. $\boldmath{\sum\limits_{i=1}^r \deg f_i = c(\mathfrak{g}) - \deg p_{\mathfrak{g}}}$ (54 cases)}\\
Then $\mathfrak{g}$ is coregular since $Y(\mathfrak{g})$ is freely generated by $f_1,\ldots, f_r$ by Theorem 14.  We have two subclasses:\\
Ia. $\mathfrak{g}$ is singular (19 cases)\\
Ib. $\mathfrak{g}$ is nonsingular (35 cases)\\
We describe the procedure by giving a nontrivial example for each case.\\
\ \\
{\bf Example 32.} (See Remark 29)\\
Let $\mathfrak{g}$ be the nilpotent Lie algebra $\mathfrak{g}_{7,2.4} \cong (23457C) \cong R_{13}$ (\# 101) with basis $\{x_1,\ldots,x_7\}$ and nonzero brackets: $[x_1, x_2] = x_3$, $[x_1, x_3] = x_4$, $[x_1, x_4] = x_5$, $[x_1, x_5] = x_6$, $[x_2, x_5] = -x_7$, $[x_3, x_4] = x_7$.\\
\ \\
Clearly, $i(\mathfrak{g}) = 3 > 2 = \rank \ \mathfrak{g}$ (the 2nd index of $\mathfrak{g}_{7,2.4}$).  Hence $c(\mathfrak{g}) = \frac{1}{2} (7 + 3) = 5$.\\
The symmetric bilinear form $b$ on $\mathfrak{g}$ with nonzero entries 
$$b (x_4, x_4) = b(x_2, x_6) = b(x_1, x_7) = 1 \ \mbox{and}\ b(x_3, x_5) = -1$$
is invariant and nondegenerate. \\
So, $\mathfrak{g}$ is quadratic and hence also quasi quadratic (i.e. $F(\mathfrak{g}) = \mathfrak{g})$.  In particular $\mathfrak{g}$ has no CP's.  Next, we consider the flag of ideals $\mathfrak{g}_0 = \{0\}$, $\mathfrak{g}_1 = \langle x_7\rangle$, $\mathfrak{g}_2 = \langle x_6, x_7\rangle,\ldots, \mathfrak{g}_6 = \langle x_2,\ldots, x_7\rangle$, $\mathfrak{g}_7 = \mathfrak{g}$ for which $[\mathfrak{g}, \mathfrak{g}_j] \subset \mathfrak{g}_{j-1}$, $j:1,\ldots, 7$.\\
$Y(\mathfrak{g})$ contains the following homogeneous elements:
$$f_1 = x_6, f_2 = x_7, f_3 = x_4^2 - 2x_3x_5 + 2x_2 x_6 + 2x_1x_7$$
($f_3$ is the Casimir element w.r.t. b).  By Theorem 31 they form algebraically independent generators of the quotient field of $Y(\mathfrak{g})$.  One verifies that $p_{\mathfrak{g}} = x_7$, hence $\mathfrak{g}$ is singular.  We also observe that
$$\sum\limits_{i=1}^3 \deg f_i = 4 = 5-1 = c(\mathfrak{g}) - \deg p_{\mathfrak{g}}$$
By Theorem 14 we may deduce that
$$Y(\mathfrak{g}) = k[f_1, f_2, f_3]$$
i.e. $\mathfrak{g}$ is coregular.  Next, we take $\xi = x_7^\ast \in \mathfrak{g}_{\reg}^\ast$.  Then the Mishchenko-Fomenko algebra
$$Y_{\xi} (\mathfrak{g}) = k[x_1, x_6, x_7, f_3]$$
has transcendence degree $4 < c(\mathfrak{g})$, as expected by Theorem 21, so it is not complete.  Moreover, it is not maximal Poisson commutative since it is strictly contained in the Poisson commutative subalgebra
$$k[x_1, x_6, x_7, x_5^2 - 2x_4x_6, f_3]$$
In order to follow the approach of (2) of Examples 19 we need to consider:\\
$Y{(\mathfrak{g}_0)} = k$, $Y{(\mathfrak{g}_1)} = k[x_7]$, $Y{(\mathfrak{g}_2)} = k[x_6, x_7] = Y(\mathfrak{g}_6)$, $Y(\mathfrak{g}_3) = k[x_5, x_6, x_7] = Y(\mathfrak{g}_5)$, $Y(\mathfrak{g}_4) = k[x_4, x_5, x_6, x_7]$, $Y(\mathfrak{g}_7) = Y(\mathfrak{g}) = k[x_6, x_7, f_3]$.  The subalgebra generated by the union of these is:
$$V(\mathfrak{g}) = k[x_4, x_5, x_6, x_7, f_3] = k[x_4, x_5, x_6, x_7, x_1x_7 + x_2x_6 - x_3x_5]$$
The Jacobian locus of its generators is:
$$\{\xi \in \mathfrak{g}^\ast \mid \xi (x_5) = \xi(x_6) = \xi(x_7) = 0\}$$
which has codimension 3.  Then combining Theorem 16 and (ii) of Proposition 15 we deduce that $M = V(\mathfrak{g})$ is a strongly complete Poisson commutative subalgebra of $S(\mathfrak{g})$.  Note that $M$ satisfies the Milovanov conjecture.  In addition $M$ is a polynomial subalgebra for which $\{\mathfrak{g}, M\} \subset M$ and $Q(M)$ is a maximal subfield of the Poisson field $R(\mathfrak{g})$ (see (2) of Examples 19).  In view of Remark 2 and (2) of Examples 19 we may conclude that:
$$Z(U(\mathfrak{g})) = k[x_6, x_7, x_4^2-2x_3x_5 + 2x_2x_6 + 2x_1x_7]$$
and that
$$s(M) = k[x_4, x_5, x_6, x_7, x_1x_7 + x_2x_6 - x_3x_5]$$
is a polynomial, maximal commutative subalgebra of $U(\mathfrak{g})$ of transcendence degree $c(\mathfrak{g})$ for which $[\mathfrak{g},s(M)] \subset s(M)$.\\
\ \\
{\bf Example 33.} (See Remark 29)\\
Let $\mathfrak{g}$ be the nilpotent Lie algebra $\mathfrak{g}_{7,1.1(iii)}\cong (23457G) \cong R_{11}$ (\# 136) with basis $\{x_1,\ldots, x_7\}$ and nonzero brackets: $[x_1, x_2] = x_3$, $[x_1, x_3] = x_4$, $[x_1, x_4] = x_5$, $[x_1, x_5] = x_6$, $[x_2, x_3] = x_5$, $[x_2, x_4] = x_6$, $[x_2, x_5] = -x_7$, $[x_3, x_4] = x_7$.\\
Clearly, $i(\mathfrak{g}) = 3 > 1 = \rank\ \mathfrak{g}$ and $c(\mathfrak{g}) = \frac{1}{2} (7+3) = 5$.\\
One verifies that $\mathfrak{g}$ is quasi quadratic (i.e. $F(\mathfrak{g}) = \mathfrak{g})$ and therefore $\mathfrak{g}$ has no CP's.  Consider the flag of ideals $\mathfrak{g}_0 = \{0\}$, $\mathfrak{g}_1 =\langle x_7\rangle$, $\mathfrak{g}_2 = \langle x_6, x_7\rangle, \ldots, \mathfrak{g}_6 = \langle x_2,\ldots, x_7\rangle$, $\mathfrak{g}_7 = \mathfrak{g}$ for which $[\mathfrak{g}, \mathfrak{g}_j] \subset \mathfrak{g}_{j-1}$, $j:1,\ldots, 7$.\\
\ \\
$Y(\mathfrak{g})$ contains the following homogeneous elements:
$$f_1 = x_6, f_2 = x_7, f_3 = 2x_5^3 - 3x_4^2 x_7 - 6x_4x_5x_6 + 6x_3x_6^2 + 6x_3x_5x_7 - 6x_2x_6x_7 - 6x_1x_7^2$$
By Theorem 31 they form algebraically independent generators of the quotient field of $Y(\mathfrak{g})$. One verifies that $\codim \ \mathfrak{g}_{\sing}^\ast = 3$.  In particular $\mathfrak{g}$ is nonsingular, so $\deg p_{\mathfrak{g}} = 0$ by Remark 5. Since
$$\sum\limits_{i=1}^3 \deg f_i = 5 = c(\mathfrak{g}) - \deg p_\mathfrak{g}$$
we obtain that $Y(\mathfrak{g}) = k[f_1, f_2, f_3]$ by Theorem 14, so $\mathfrak{g}$ is coregular.  Using the same reasoning as in the previous example, we obtain that
\begin{eqnarray*}
M &=& V(\mathfrak{g}) = k[x_4, x_5, x_6, x_7, f_3]\\
&=& k [x_4, x_5, x_6, x_7, x_3x_6^2 + x_3 x_5x_7 - x_2x_6x_7 - x_1x_7^2]
\end{eqnarray*}
is a polynomial, strongly complete Poisson commutative subalgebra of $S(\mathfrak{g})$, for which $\{\mathfrak{g},M\} \subset M$ and $Q(M)$ is a maximal subfield of the Poisson field $R(\mathfrak{g})$.  These results can be literally carried over to the enveloping algebra $U(\mathfrak{g})$.  Indeed, each of the monomials appearing in $f_3$ is a product of commuting variables.  By Remark 2 and (2) of Examples 19 we have that 
$$Z(U(\mathfrak{g})) = k[x_6, x_7, f_3]\ \mbox{and}\ s(M) = k[x_4, x_5, x_6, x_7, f_3]$$
where $s(M)$  has similar properties as $M$.\\
The same phenomenon occurs in all other Lie algebras of our list.\\
Next, we take $\xi = x_5^\ast \in \mathfrak{g}_{\reg}^\ast$.  The $\xi$-shifts of $f_3$ are $f_3, f_3^1 = 6(x_5^2 - x_4x_6 + x_3x_7)$ and $f_3^2 = 6x_5$.  Then, by Theorem 24, $Y_{\xi} (\mathfrak{g}) = k[x_6, x_7, f_3, f_3^1, f_3^2]\\
= k[x_5, x_6, x_7, x_4x_6 - x_3x_7, x_4^2 x_7 - 2x_3x_6^2 + 2x_2x_6x_7 + 2x_1x_7^2]$\\
 is a polynomial, strongly complete Poisson commutative subalgebra of $S(\mathfrak{g})$.\\
 However, $\{\mathfrak{g}, Y_\xi (\mathfrak{g})\} \not\subset Y_\xi(\mathfrak{g})$ (indeed $\{x_1, x_3x_7 - x_4x_6\} = x_4x_7 - x_5x_6 \notin Y_\xi(\mathfrak{g})$).\\
 Finally, following a suggestion by Oksana Yakimova, we decompose $f_3$ as follows: 
 $$f_3 = 2x_5^3 + 6x_6u - 3x_7v$$
where $u = x_3x_6 - x_4x_5$ and $v = x_4^2 + 2x_1x_7 + 2x_2x_6 - 2x_3x_5$.\\
 It turns out (by Theorem 16 and (ii) of Proposition 15) that
 $$M_1 = k[x_5, x_6, x_7, u, v]$$
 is a polynomial, strongly complete Poisson commutative subalgebra of $S(\mathfrak{g})$, of which the generators have degrees at most two.\\
\ \\
 {\bf II. ${\boldmath\sum\limits_{i=1}^r \deg f_i > c(\mathfrak{g})}$ (23 cases)}\\
 It turns out that in each case $\mathfrak{g}$ is not coregular.  The minimal generators of $Y(\mathfrak{g})$ satisfy only one relation, except for $\mathfrak{g}_{7,3.2}$ (\# 158) (63 relations) and $\mathfrak{g}_{7,2.3}$ (\# 159) (168 relations).\\
\ \\
{\bf Example 34.} Let $\mathfrak{g}$ be the nilpotent Lie algebra $\mathfrak{g}_{7,1.01(i)} \cong (12357B) \cong R_{43}$ (\# 154) with basis $\{x_1, \ldots, x_7\}$ and nonzero brackets $[x_1, x_2] = x_4$, $[x_1, x_4] = x_5$, $[x_1, x_5] = x_6$, $[x_1, x_6] = x_7$, $[x_2, x_3] = x_5 + x_7$, $[x_3, x_4] = -x_6$, $[x_3, x_5] = -x_7$.\\
It is easy to verify that $i(\mathfrak{g}) = 3 > 1 = \rank\ \mathfrak{g}$ so $c(\mathfrak{g}) = \frac{1}{2} (7 + 3) = 5$ and also $\mathfrak{h} = \langle x_2, x_4, x_5, x_6, x_7\rangle$ is a CP of $\mathfrak{g}$.  On the other hand, $\codim \ \mathfrak{g}_{\sing}^\ast = 3$.  Consequently, $\mathfrak{g}$ is not coregular by (3) of Theorem 13.  In order to find the generators of $Y(\mathfrak{g})$ we use Dixmier's methods of [D2, pp. 328-329].  So, consider the flag of ideals $\mathfrak{g}_0 = \{0\}$, $\mathfrak{g}_1 = kx_7$, $\mathfrak{g}_2 = \langle x_6, x_7\rangle, \ldots, \mathfrak{g}_6 = \langle x_2, \ldots, x_7\rangle$, $\mathfrak{g}_7 = \mathfrak{g}$, for which $[\mathfrak{g},\mathfrak{g}_i] \subset \mathfrak{g}_{i-1}$, $i : 1,\ldots, 7$.\\
One verifies that $Y(\mathfrak{g})$ contains the following $x_7, f = x_6^3 - 3x_5x_6x_7 + 3x_4x_7^2$,\\
$g = x_6^4 - 4x_5x_6^2 x_7 - 2x_6^2x_7^2 + 2x_5^2x_7^2 + 4x_4 x_6 x_7^2 + 4x_5x_7^3 - 4x_2x_7^3$. \\
They satisfy the conditions of Theorem 31 and therefore form algebraically independent (over $k$) generators of the quotient field of $Y(\mathfrak{g})$ and
$$Y(\mathfrak{g}) \subset k[x_7, f, g, x_7^{-1}]$$
Next, put $h = (f^4 - g^3 - 6x_7^2 f^2g)/x_7^3 \in Y(\mathfrak{g})$.  We claim that $Y(\mathfrak{g}) = k[x_7, f, g, h]$.  For this we will need the following lemmas.\\
\ \\
{\bf Lemma A.} Let $P$ be a polynomial in 3 variables $X$, $Y$, $Z$ with coefficients in $k$.  If $P(f,g,h)$ is divisible by $x_7$ in $S(\mathfrak{g})$, then $P(X,Y,Z)$ is divisible by $X^4 - Y^3$ in $k [X,Y,Z]$.\\
\ \\
{\bf Proof.} Let $I$ be the ideal of $S(\mathfrak{g})$ generated by $x_7$.  We identify the quotient $S(\mathfrak{g})/I$ with $k [x_1,\ldots, x_6]$. Then the canonical images of $f$, $g$, $h$ are
$$x_6^3, x_6^4, 4 (x_5^3 - 3x_4 x_5 x_6 + 3x_2 x_6^2) x_6^6$$
By§ assumption we obtain
$$P(x_6^3, x_6^4, 4(x_5^3 - 3x_4x_5x_6 + 3x_2x_6^2) x_6^6) = 0$$
We now decompose $P$ as follows:
$$P(X,Y,Z) = P_r(X,Y) Z^r + P_{r-1} (X,Y) Z^{r-1} + \ldots + P_0 (X,Y)$$
Considering the terms in $x_5$, we see that for all $i : P_i (x_6^3, x_6^4) = 0$.\\
So each $P_i$ is divisible by $X^4 - Y^3$.  Consequently, the same holds for $P$.\\
\ \\
{\bf Lemma B.} Let $q \in S(\mathfrak{g})$ such that $x_7q \in k[x_7, f,g,h]$.  Then $q \in k[x_7, f, g, h]$.\\
\ \\
{\bf Proof.} By assumption there are $g_i \in k[f,g,h]$ such that
$$x_7q = x_7^rg_r + \ldots + x_7^3 g_3 + x_7^2 g_2 + x_7g_1 + g_0$$
Clearly, $g_0$ is divisible by $x_7$ in $S(\mathfrak{g})$ and hence by the previous lemma also by $f^4 - g^3$, i.e. there is a $g_0' \in k[f,g,h]$ such that
$$g_0 = (f^4 - g^3) g'_0 = (x_7^3 h + 6x_7^2 f^2 g) g'_0$$
Therefore
$$x_7q = x_7^rg_r + \ldots + x_7^3(g_3 + hg'_0) + x_7^2 (g_2 + 6f^2gg'_0) + x_7g_1$$
and finally
$$q = x_7^{r-1} g_r + \ldots + x_7^2 (g_3 + hg'_0) + x_7 (g_2 + 6f^2gg'_0) + g_1$$
which belongs to $k[x_7, f, g, h]$.\\
\ \\
We can now establish the claim.  Indeed, take $q \in Y(\mathfrak{g})$. \\
In particular, $q \in k[x_7, f, g, x_7^{-1}]$, i.e. for some $t\ x_7^t q \in k[x_7, f, g] \subset k[x_7, f, g, h]$.\\
\ \\
By applying the previous lemma $t$ times, we arrive at $q \in k[x_7, f, g, h]$.  Hence, $Y(\mathfrak{g}) \subset k[x_7, f, g, h]$.  The other inclusion is obvious.  Also, $Z(U(\mathfrak{g})) = k[x_7, f, g, h]$. \\
Finally, since $\mathfrak{h} = \langle x_2, x_4, x_5, x_6, x_7\rangle$ is a CPI of $\mathfrak{g}$ it follows that $M = k[x_2, x_4, x_5, x_6, x_7]$ is a polynomial, strongly complete Poisson commutative subalgebra of $S(\mathfrak{g})$ (see (1) of Examples 19) with generators of degree one and $\{\mathfrak{g},M\} \subset M$.\\
\ \\
{\bf Example 35.} Let $\mathfrak{g}$ be the nilpotent Lie algebra $\mathfrak{g}_{7, 1.21} \cong (12457F) \cong R_{36}$ (\# 150) with basis $\{x_1,\ldots, x_7\}$ and nonzero brackets:\\
$[x_1, x_2] = x_4$, $[x_1, x_4] = x_5$, $[x_1, x_5] = x_6$, $[x_2, x_3] = x_6$, $[x_2, x_4] = x_6$, $[x_2, x_6] = x_7$, $[x_4, x_5] = -x_7$.\\
$i(\mathfrak{g}) = 3 > 1 = \rank\ \mathfrak{g}$ and $c(\mathfrak{g}) = \frac{1}{2} (7 + 3) = 5$.\\
One verifies that $F(\mathfrak{g}) = \langle x_1, x_3, x_4, x_5, x_6, x_7\rangle$, which is not commutative.  Therefore $\mathfrak{g}$ has no CP's.  Note that $Z(\mathfrak{g}) = \langle x_7\rangle$ and that $\codim \ \mathfrak{g}_{\sing}^\ast = 2$.  Hence $\deg p_{\mathfrak{g}} = 0$ by Remark 5. Consider the flag of ideals $\mathfrak{g}_0 = \{0\}$, $\mathfrak{g}_1 = \{x_7\}$, $\mathfrak{g}_2 = \langle x_6, x_7\rangle, \ldots,\\
\mathfrak{g}_6 = \langle x_2,\ldots, x_7\rangle$, $\mathfrak{g}_7 = \mathfrak{g}$, for which $[\mathfrak{g}, \mathfrak{g}_j] \subset \mathfrak{g}_{j-1}$, $j : 1,\ldots, 7$.\\
$Y(\mathfrak{g})$ contains the following homogeneous elements:
$$x_7, f = x_6^2 - 2x_3 x_7, g = 2x_6^3 - 6x_4x_6x_7 + 3x_5^2x_7 - 6x_1x_7^2$$
By Theorem 31 they form algebraically independent generators of the quotient field of $Y(\mathfrak{g})$.  It is easy to check that any element of $Y(\mathfrak{g})$ of degree at most 2 is a linear combination of $x_7$, $x_7^2$ and $f$.  Now suppose $\mathfrak{g}$ were to be coregular.  By (1) of Theorem 13 the equality
$$3i(\mathfrak{g}) + 2\deg p_{\mathfrak{g}} = 9 = \dim \mathfrak{g} + 2\dim Z(\mathfrak{g})$$
would imply that $Y(\mathfrak{g}) = k[x_7, f]$, which is impossible.  So, $\mathfrak{g}$ is not coregular.\\
Next, put $h = (4f^3 - g^2)/x_7 \in Y(\mathfrak{g})$.\\
Using a similar argument as in the previous example we obtain that
$$Y(\mathfrak{g}) = k[x_7, f, g, h]$$
Next, we observe that $Y(\mathfrak{g}_0) = k$, $Y(\mathfrak{g}_1) = k[x_7], Y(\mathfrak{g}_2) = k[x_6, x_7] = Y(\mathfrak{g}_4)$,\\
$Y(\mathfrak{g}_3) = k[x_5, x_6, x_7], Y(\mathfrak{g}_5) = k[x_3, x_6, x_7]$, $Y(\mathfrak{g}_6) = k[x_7, f], \\
Y(\mathfrak{g}_7) = Y(\mathfrak{g}) = k [x_7, f, g, h]$.\\
The subalgebra generated by the union of all these is:
\begin{eqnarray*}
V(\mathfrak{g}) &=& k[x_3, x_5, x_6, x_7, f, g, h]\\
&=& k[x_3, x_5, x_6, x_7, (x_1x_7 + x_4x_6) x_7, h]
\end{eqnarray*}
Put $M = k[x_3, x_5, x_6, x_7, x_1x_7 + x_4x_6] \subset S(\mathfrak{g})$.\\
The Jacobian locus of its generators is
$$\{\xi \in \mathfrak{g}^\ast \mid \xi(x_6) = 0 = \xi(x_7)\}$$
which has codimension 2.  Combining Theorem 16 and (ii) of Proposition 15 we see that $M$ is a polynomial, strongly complete Poisson commutative subalgebra of $S(\mathfrak{g})$.  Clearly, $f, g \in M$.  Since $x_7$ and $hx_7 = 4f^3 - g^2 \in M$ we obtain that $h \in M$ by (i) of Remark 17 because $M$ is a maximal, Poisson commutative subalgebra of $S(\mathfrak{g})$.  Consequently, $V(\mathfrak{g}) \subset M$.  On the other hand, $M \subset V(\mathfrak{g})_{x_7} \cap S(\mathfrak{g}) \subset Q(V(\mathfrak{g})) \cap S(\mathfrak{g})$.\\
Hence, $M = V(\mathfrak{g})_{x_7} \cap S(\mathfrak{g}) = Q(V(\mathfrak{g})) \cap S(\mathfrak{g}) = V(\mathfrak{g})'$ by the maximality of $M$.  This means that we are in the situation of (2) of Examples 19.  In particular, $\{g,M\} \subset M$, $Z(U(\mathfrak{g})) = k [x_7, f, g, h]$ and $s(M)$ is a polynomial, maximal commutative subalgebra of $U(\mathfrak{g})$.  Finally, we observe that the Milovanov conjecture is satisfied and also that $M \subset S(F(\mathfrak{g}))$.  So, $M$ is also a strongly complete Poisson commutative subalgebra of $S(F(\mathfrak{g}))$.\\
\ \\
{\bf 3.2. The list of 7-dimensional indecomposable nilpotent Lie algebras with $i(\mathfrak{g}) > \rank\ \mathfrak{g}$}\\
The main purpose is to describe for each Lie algebra $\mathfrak{g}$ the Poisson center $Y = Y(\mathfrak{g})$ and also to give a polynomial, strongly complete Poisson commutative subalgebra $M \subset S(\mathfrak{g})$ for which $\{\mathfrak{g},M\} \subset M$ and the quotient field $Q(M)$ is a maximal Poisson commutative subfield of the Poisson field $R(\mathfrak{g})$.  As pointed out earlier (see Example 33) these results may also be interpreted as taking place in the enveloping algebra $U(\mathfrak{g})$.  So $Y$ can be regarded as the center $Z(U(\mathfrak{g}))$ and $M$ as a maximal commutative subalgebra of $U(\mathfrak{g})$ with similar properties.  Moreover, the generators of $M$ have degrees at most two, except for the Lie algebras 134, 135 and 136.  For the latter we offer an alternative polynomial, strongly complete subalgebra $M_1 \subset S(\mathfrak{g})$ with generators of degree at most two.  Furthermore, we also give the Frobenius semi-radical $F = F(\mathfrak{g})$ and if it exists a CP-ideal (CPI) of $\mathfrak{g}$.  Other abbreviations are: $i = i(\mathfrak{g})$, $r = \rank \ \mathfrak{g}$, $c = c(\mathfrak{g})$, $p = p_{\mathfrak{g}}$ the fundamental semi-invariant, $cod = \codim\ (\mathfrak{g}^{\ast}_{\sing})$, SQ.I. = square integrable, $Y_{\xi}$ is the Mishchenko-Fomenko algebra.\\
\ \\
{\bf Notation.} We use the same notation as Magnin [Ma2, Ma3], but we will also include the notation of Seeley ((...)) [Se] and Romdhani ($R_-$) [R].  Since this is the continuation of [O5,5], where 82 different cases were listed, we start numbering at 83.\\
\newpage

{\bf I. $\mathfrak{g}$ is coregular} (54 cases) \\
\ \\
{\bf Ia. $\mathfrak{g}$ is singular} (i.e. $\codim \ \mathfrak{g}_{\sing}^\ast = 1$, 19 cases)
\begin{itemize}
\item[83.] $\mathfrak{g}_{7,0.1} \cong (123457F) \cong R_2$\\
$[x_1,x_2] = x_3$, $[x_1, x_3] = x_4$, $[x_1, x_4] = x_5$, $[x_1,x_5] = x_6$, $[x_1, x_6] = x_7$,\\
$[x_2, x_3] = x_6$, $[x_2, x_4] = x_7$, $[x_2, x_5] = x_7$, $[x_3, x_4] = -x_7$.\\
SQ.I., $i = 1$, $r= 0$, $c=4$, $F = \langle x_7\rangle$, $CPI = \langle x_4, x_5, x_6, x_7\rangle$,\\
$p = x_7^3$, $Y = k[x_7]$, $M = k[x_4, x_5, x_6, x_7]$.

\item[84.] $\mathfrak{g}_{7,0.4(\lambda)} \cong (12457N)(\xi \neq 1) \cong R_{14}^\lambda$\\
$[x_1,x_2] = x_3$, $[x_1, x_3] = x_4$, $[x_1, x_4] = \lambda x_7 + x_6$, $[x_1,x_5] = x_7$, $[x_1, x_6] = x_7$,\\
$[x_2, x_3] = x_5$, $[x_2, x_4] = x_7$, $[x_2, x_5] = x_6$, $[x_3, x_5] = x_7$.\\
SQ.I., $i = 1$, $r= 0$, $c=4$, $F = \langle x_7\rangle$, $CPI = \langle x_4, x_5, x_6, x_7\rangle$,\\
$p = x_7^3$, $Y = k[x_7]$, $M = k[x_4, x_5, x_6, x_7]$.

\item[85.] $\mathfrak{g}_{7,0.6} \cong (12457J) \cong R_{20}$\\
$[x_1,x_2] = x_3$, $[x_1, x_3] = x_4$, $[x_1, x_4] = x_7$, $[x_1,x_5] = x_6$, $[x_1, x_6] = x_7$,\\
$[x_2,x_3] = x_5$, $[x_2, x_4] = x_6$, $[x_2, x_5] = x_7$, $[x_3,x_4] = x_7$.\\
SQ.I., $i = 1$, $r= 0$, $c=4$, $F = \langle x_7\rangle$, $CPI = \langle x_4, x_5, x_6, x_7\rangle$,\\
$p = x_7^3$, $Y = k[x_7]$, $M = k[x_4, x_5, x_6, x_7]$.

\item[86.] $\mathfrak{g}_{7,0.7} \cong (13457I) \cong R_{25}$\\
$[x_1,x_2] = x_3$, $[x_1, x_3] = x_4$, $[x_1, x_4] = x_7$, $[x_1,x_5] = x_7$, $[x_1, x_6] = x_7$,\\
$[x_2, x_3] = x_5$, $[x_2, x_4] = x_7$, $[x_2, x_5] = x_6$, $[x_3, x_5] = x_7$.\\
SQ.I., $i = 1$, $r= 0$, $c=4$, $F = \langle x_7\rangle$, $CPI = \langle x_4, x_5, x_6, x_7\rangle$,\\
$p = x_7^3$, $Y = k[x_7]$, $M = k[x_4, x_5, x_6, x_7]$.

\item[87.] $\mathfrak{g}_{7,0.8} \cong (12457G) \cong R_{34}$\\
$[x_1,x_2] = x_4$, $[x_1, x_3] = x_7$, $[x_1, x_4] = x_5$, $[x_1,x_5] = x_6$,\\
$[x_2, x_3] = x_6$, $[x_2, x_4] = x_6$, $[x_2, x_6] = x_7$, $[x_4, x_5] = -x_7$.\\
SQ.I., $i = 1$, $r= 0$, $c=4$, $F = \langle x_7\rangle$, $CPI = \langle x_3, x_5, x_6, x_7\rangle$,\\
$p = x_7^3$, $Y = k[x_7]$, $M = k[x_3, x_5, x_6, x_7]$.

\item[88.] $\mathfrak{g}_{7,2.40} \cong (357C)$\\
$[x_1,x_2] = x_3$, $[x_1, x_3] = x_5$, $[x_1, x_4] = x_7$, $[x_2,x_3] = x_6$, $[x_2, x_4] = x_5$.\\
SQ.I., $i = 3$, $r= 2$, $c=5$, $F = \langle x_5, x_6, x_7\rangle$, $CPI = \langle x_3, x_4, x_5, x_6, x_7\rangle$,\\
$p = x_5^2- x_6x_7$, $Y = k[x_5, x_6, x_7]$, $M = k[x_3, x_4, x_5, x_6, x_7]$.

\item[89.] $\mathfrak{g}_{7,2.19} \cong (2457G) \cong R_{61}$\\
$[x_1,x_2] = x_4$, $[x_1, x_3] = x_6$, $[x_1, x_4] = x_5$, $[x_1,x_5] = x_7$, $[x_2, x_4] = x_6$.\\
$i = 3$, $r= 2$, $c=5$, $F = \langle x_3, x_5, x_6, x_7\rangle$, $CPI = \langle x_3, x_4, x_5, x_6, x_7\rangle$,\\
$p = x_6$, $Y = k[x_6, x_7, x_5x_6 - x_3x_7]$, $M = k[x_3, x_4, x_5, x_6, x_7]$.

\item[90.] $\mathfrak{g}_{7,2.20} \cong (2457F) \cong R_{76}$\\
$[x_1,x_2] = x_4$, $[x_1, x_3] = x_5$, $[x_1, x_5] = x_6$, $[x_1,x_6] = x_7$, $[x_3, x_5] = x_7$.\\
$i = 3$, $r= 2$, $c=5$, $F = \langle x_2, x_4, x_6, x_7\rangle$, $CPI = \langle x_2, x_4, x_5, x_6, x_7\rangle$,\\
$p = x_7$, $Y = k[x_4, x_7, x_4x_6 - x_2x_7]$, $M = k[x_2, x_4, x_5, x_6, x_7]$.

\item[91.] $\mathfrak{g}_{7,2.39} \cong (247L)$\\
$[x_1,x_2] = x_4$, $[x_1, x_3] = x_5$, $[x_1, x_4] = x_6$, $[x_1,x_5] = x_7$, $[x_2, x_3] = x_6$.\\
$i = 3$, $r= 2$, $c=5$, $F = \langle x_4, x_5, x_6, x_7\rangle$, $CPI = \langle x_3, x_4, x_5, x_6, x_7\rangle$,\\
$p = x_6$, $Y = k[x_6, x_7, x_5x_6 - x_4x_7]$, $M = k[x_3, x_4, x_5, x_6, x_7]$.

\item[92.] $\mathfrak{g}_{7,2.27} \cong (257I) \cong R_{106}$\\
$[x_1,x_2] = x_5$, $[x_1, x_3] = x_7$, $[x_1, x_5] = x_6$, $[x_2,x_4] = x_7$, $[x_2, x_5] = x_7$.\\
$i = 3$, $r= 2$, $c=5$, $F = \langle x_3, x_4 - x_5, x_6, x_7\rangle$, $CPI = \langle x_3, x_4, x_5, x_6, x_7\rangle$,\\
$p = x_7$, $Y = k[x_6, x_7, x_3x_6 + x_4x_7 - x_5x_7]$, $M = k[x_3, x_4, x_5, x_6, x_7]$.

\item[93.] $\mathfrak{g}_{7,1.16} \cong (2457D)$\\
$[x_1,x_2] = x_3$, $[x_1, x_3] = x_5$, $[x_1, x_4] = x_6$, $[x_1,x_5] = x_7$, $[x_1, x_6] = x_7$,\\
$[x_2, x_3] = x_7$.\\
$i = 3$, $r= 1$, $c=5$, $F = \langle x_4, x_5, x_6, x_7\rangle$, $CPI = \langle x_3, x_4, x_5, x_6, x_7\rangle$,\\
$p = x_7$, $Y = k[x_7, x_5-x_6, x_6^2 - 2x_4x_7]$, $M = k[x_3, x_4, x_5, x_6, x_7]$.

\item[94.] $\mathfrak{g}_{7,2.7} \cong (23457A) \cong R_{31}$\\
$[x_1,x_2] = x_3$, $[x_1, x_3] = x_4$, $[x_1, x_4] = x_6$, $[x_1,x_6] = x_7$, $[x_2, x_3] = x_5$.\\
$i = 3$, $r= 2$, $c=5$, $F = \langle x_4, x_5, x_6, x_7\rangle$, $CPI = \langle x_3, x_4, x_5, x_6, x_7\rangle$,\\
$p = x_5$, $Y = k[x_5, x_7, x_6^2 - 2x_4x_7]$, $M = k[x_3, x_4, x_5, x_6, x_7]$.

\item[95.] $\mathfrak{g}_{7,2.21} \cong (2457C) \cong R_{77}$\\
$[x_1,x_2] = x_4$, $[x_1, x_3] = x_5$, $[x_1, x_4] = x_6$, $[x_1,x_6] = x_7$, $[x_2, x_3] = x_7$.\\
$i = 3$, $r= 2$, $c=5$, $F = \langle x_4, x_5, x_6, x_7\rangle$, $CPI = \langle x_3, x_4, x_5, x_6, x_7\rangle$,\\
$p = x_7$, $Y = k[x_5, x_7, x_6^2 - 2x_4x_7]$, $M = k[x_3, x_4, x_5, x_6, x_7]$.

\item[96.] $\mathfrak{g}_{7,2.43} \cong (247C)$\\
$[x_1,x_2] = x_4$, $[x_1, x_3] = x_5$, $[x_1, x_4] = x_7$, $[x_1,x_5] = x_6$, $[x_3, x_5] = x_7$.\\
$i = 3$, $r= 2$, $c=5$, $F = \langle x_2, x_4, x_6, x_7\rangle$, $CPI = \langle x_2, x_4, x_5, x_6, x_7\rangle$,\\
$p = x_7$, $Y = k[x_6, x_7, x_4^2 - 2x_2x_7]$, $M = k[x_2, x_4, x_5, x_6, x_7]$.

\item[97.] $\mathfrak{g}_{7,2.11} \cong (2457E) \cong R_{58}$\\
$[x_1,x_2] = x_4$, $[x_1, x_4] = x_5$, $[x_1, x_5] = x_7$, $[x_2,x_3] = x_6$, $[x_2, x_4] = x_6$.\\
$i = 3$, $r= 2$, $c=5$, $F = \langle x_3-x_4, x_5, x_6, x_7\rangle$, $CPI = \langle x_3, x_4, x_5, x_6, x_7\rangle$,\\
$p = x_6$, $Y = k[x_6, x_7, x_5^2 + 2x_3x_7 - 2x_4x_7]$, $M = k[x_3, x_4, x_5, x_6, x_7]$.

\item[98.] $\mathfrak{g}_{7,2.18} \cong (2457H) \cong R_{60}$\\
$[x_1,x_2] = x_4$, $[x_1, x_4] = x_5$, $[x_1, x_5] = x_7$, $[x_2,x_3] = x_7$, $[x_2, x_4] = x_6$.\\
$i = 3$, $r= 2$, $c=5$, $F = \langle x_3, x_4, x_5, x_6, x_7\rangle = CPI$,\\
$p = x_7$, $Y = k[x_6, x_7, x_5^2 + 2x_3x_6 - 2x_4x_7]$, $M = k[x_3, x_4, x_5, x_6, x_7]$.

\item[99.] $\mathfrak{g}_{7,2.44} \cong (247N)$\\
$[x_1,x_2] = x_4$, $[x_1, x_3] = x_5$, $[x_1, x_5] = x_7$, $[x_2,x_3] = x_6$, $[x_2, x_4] = x_7$.\\
$i = 3$, $r= 2$, $c=5$, $F = \langle x_3, x_4, x_5, x_6, x_7\rangle = CPI$,\\
$p = x_7$, $Y = k[x_6, x_7, x_5^2 + 2x_4x_6 - 2x_3x_7]$, $M = k[x_3, x_4, x_5, x_6, x_7]$.

\item[100.] $\mathfrak{g}_{7,2.12} \cong (247E) \cong R_{91}$\\
$[x_1,x_2] = x_4$, $[x_1, x_3] = x_5$, $[x_1, x_4] = x_6$, $[x_2,x_4] = x_7$, $[x_3, x_5] = x_7$.\\
$i = 3$, $r= 2$, $c=5$, $F = \langle x_1, x_2, x_4, x_5, x_6, x_7\rangle$;, no $CP's$,\\
$p = x_7$, $Y = k[x_6, x_7, x_4^2 + x_5^2 - 2x_2x_6 + 2x_1 x_7]$, $M = k[x_4, x_5, x_6, x_7, x_2x_6 - x_1x_7]$.

\item[101.] $\mathfrak{g}_{7,2.4} \cong (23457C) \cong R_{13}$ (quadratic, see Example 32)\\
$[x_1,x_2] = x_3$, $[x_1, x_3] = x_4$, $[x_1, x_4] = x_5$, $[x_1,x_5] = x_6$, $[x_2, x_5] = -x_7$,\\
$[x_3, x_4] = x_7$.\\
$i = 3$, $r= 2$, $c=5$, $F = \mathfrak{g}_{7,2.4}$, no $CP's$, $p = x_7$, $Y = k[x_6, x_7, f]$,\\
$M = k[x_4, x_5, x_6, x_7, f]$ where $f = x_4^2 - 2x_3x_5 + 2x_2x_6 + 2x_1x_7$.
\end{itemize}

{\bf Ib. $\mathfrak{g}$ is nonsingular} (35 cases)\\
\ \\
$i(\mathfrak{g}) = 3$; $\codim \ \mathfrak{g}_{sing}^\ast = 2$, except 3 for $\mathfrak{g}_{7,1.1(iii)} (\# 136)$.
\begin{itemize}
\item[102.] $\mathfrak{g}_{7,1.2(i_\lambda), \lambda = 1} \cong (1357Q) \cong R_{52}^1$\\
$[x_1,x_2] = x_4$, $[x_1, x_3] = x_6$, $[x_1, x_4] = x_5$, $[x_1,x_5] = x_7$, $[x_2, x_3] = x_5$, \\
$[x_2, x_4] = x_6$, $[x_2, x_6] = x_7$.\\
$i = 3$, $r= 1$, $c=5$, $F = \langle x_3, x_4, x_5, x_6, x_7\rangle = CPI$,\\
$Y = k[x_7, x_5x_6 - x_3x_7, x_5^2 + x_6^2 - 2x_4x_7]$, $M = k[x_3, x_4, x_5, x_6, x_7]$.

\item[103.] $\mathfrak{g}_{7,2.1(i_\lambda), \lambda = 1} \cong (1357M)(\xi =1) \cong R_{64}^1$\\
$[x_1,x_2] = x_4$, $[x_1, x_3] = x_5$, $[x_1, x_4] = x_6$, $[x_1,x_6] = x_7$, $[x_2, x_3] = x_6$,\\
$[x_2, x_5] = x_7$.\\
$i = 3$, $r= 2$, $c=5$, $F = \langle x_3, x_4, x_5, x_6, x_7\rangle = CPI$,\\
$Y = k[x_7, x_5x_6 - x_3x_7, x_6^2 - 2x_4x_7]$, $M = k[x_3, x_4, x_5, x_6, x_7]$.

\item[104.] $\mathfrak{g}_{7,2.41} \cong (1357O)$\\
$[x_1,x_2] = x_3$, $[x_1, x_3] = x_5$, $[x_1, x_4] = x_6$, $[x_1,x_6] = x_7$, $[x_2, x_3] = x_6$,\\
$[x_2, x_5] = x_7$.\\
$i = 3$, $r= 2$, $c=5$, $F = \langle x_3, x_4, x_5, x_6, x_7\rangle = CPI$,\\
$Y = k[x_7, x_5x_6 - x_3x_7, x_6^2 - 2x_4x_7]$, $M = k[x_3, x_4, x_5, x_6, x_7]$.

\item[105.] $\mathfrak{g}_{7,2.31} \cong (1357G)\cong R_{74}$\\
$[x_1,x_2] = x_4$, $[x_1, x_4] = x_6$, $[x_1, x_6] = x_7$, $[x_2,x_3] = x_5$, $[x_2, x_5] = x_7$.\\
$i = 3$, $r= 2$, $c=5$, $F = \langle x_3, x_4, x_5, x_6, x_7\rangle = CPI$,\\
$Y = k[x_7, x_5^2 - 2x_3x_7, x_6^2 - 2x_4x_7]$, $M = k[x_3, x_4, x_5, x_6, x_7]$.

\item[106.] $\mathfrak{g}_{7,1.4} \cong (123457D)\cong R_{8}$\\
$[x_1,x_2] = x_3$, $[x_1, x_3] = x_4$, $[x_1, x_4] = x_5$, $[x_1,x_5] = x_6$, $[x_1, x_6] = x_7$,\\
$[x_2, x_3] = x_6$, $[x_2, x_4] = x_7$.\\
$i = 3$, $r= 1$, $c=5$, $F = \langle x_3, x_4, x_5, x_6, x_7\rangle = CPI$,\\
$Y = k[x_7, x_5^2 - 2x_4x_6 + 2x_3x_7, x_6^2 - 2x_5x_7]$, $M = k[x_3, x_4, x_5, x_6, x_7]$.

\item[107.] $\mathfrak{g}_{7,2.15} \cong (12457A)\cong R_{48}$\\
$[x_1,x_2] = x_4$, $[x_1, x_4] = x_5$, $[x_1, x_5] = x_6$, $[x_1,x_6] = x_7$, $[x_2, x_3] = x_6$,\\
$[x_3,x_4] = -x_7$.\\
$i = 3$, $r= 2$, $c=5$, $F = \langle x_2, x_4, x_5, x_6, x_7\rangle = CPI$,\\
$Y = k[x_7, x_5^2 - 2x_4x_6 + 2x_2x_7, x_6^2 - 2x_5x_7]$, $M = k[x_2, x_4, x_5, x_6, x_7]$.

\item[108.] $\mathfrak{g}_{7,2.38} \cong (257J)$\\
$[x_1,x_2] = x_3$, $[x_1, x_3] = x_6$, $[x_1, x_5] = x_7$, $[x_2,x_3] = x_7$, $[x_2, x_4] = x_6$.\\
$i = 3$, $r= 2$, $c=5$, $F = \langle x_3, x_4, x_5, x_6, x_7\rangle = CPI$,\\
$Y = k[x_6, x_7, x_5x_6^2 - x_3x_6x_7 + x_4x_7^2]$, $M = k[x_3, x_4, x_5, x_6, x_7]$.

\item[109.] $\mathfrak{g}_{7,2.45} \cong (257D)$\\
$[x_1,x_2] = x_5$, $[x_1, x_4] = x_7$, $[x_1, x_5] = x_6$, $[x_2,x_3] = x_7$, $[x_2, x_4] = x_6$.\\
$i = 3$, $r= 2$, $c=5$, $F = \langle x_3, x_4, x_5, x_6, x_7\rangle = CPI$,\\
$Y = k[x_6, x_7, x_3x_6^2 - x_4x_6x_7 + x_5x_7^2]$, $M = k[x_3, x_4, x_5, x_6, x_7]$.

\item[110.] $\mathfrak{g}_{7,2.22} \cong (2457I)\cong R_{79}$\\
$[x_1,x_3] = x_4$, $[x_1, x_4] = x_6$, $[x_1, x_6] = x_7$, $[x_2,x_3] = x_5$, $[x_3, x_4] = x_7$.\\
$i = 3$, $r= 2$, $c=5$, $F = \langle x_2, x_4, x_5, x_6, x_7\rangle = CPI$,\\
$Y = k[x_5, x_7, x_5x_6^2 - 2x_4x_5x_7 - 2x_2x_7^2]$, $M = k[x_2, x_4, x_5, x_6, x_7]$.

\item[111.] $\mathfrak{g}_{7,2.42} \cong (247M)$\\
$[x_1,x_2] = x_4$, $[x_1, x_3] = x_5$, $[x_1, x_5] = x_6$, $[x_2,x_3] = x_6$, $[x_2, x_4] = x_7$.\\
$i = 3$, $r= 2$, $c=5$, $F = \langle x_3, x_4, x_5, x_6, x_7\rangle = CPI$,\\
$Y = k[x_6, x_7, x_5^2x_7 - 2x_3x_6x_7 + 2x_4x_6^2]$, $M = k[x_3, x_4, x_5, x_6, x_7]$.

\item[112.] $\mathfrak{g}_{7,2.8} \cong (2457M)\cong R_{32}$\\
$[x_1,x_2] = x_3$, $[x_1, x_3] = x_4$, $[x_1, x_4] = x_6$, $[x_1,x_5] = x_7$, $[x_2, x_3] = x_5$,\\
$[x_2,x_4] = x_7$.\\
$i = 3$, $r= 2$, $c=5$, $F = \langle x_3, x_4, x_5, x_6, x_7\rangle = CPI$,\\
$Y = k[x_6, x_7, x_5^2x_6 - 2x_4x_5x_7 + 2x_3x_7^2]$, $M = k[x_3, x_4, x_5, x_6, x_7]$.

\item[113.] $\mathfrak{g}_{7,2.9} \cong (2457L)\cong R_{33}$\\
$[x_1,x_2] = x_3$, $[x_1, x_3] = x_4$, $[x_1, x_4] = x_6$, $[x_2,x_3] = x_5$, $[x_2, x_5] = x_7$.\\
$i = 3$, $r= 2$, $c=5$, $F = \langle x_3, x_4, x_5, x_6, x_7\rangle = CPI$,\\
$Y = k[x_6, x_7, x_5^2x_6 - 2x_3x_6x_7 + x_4^2x_7]$, $M = k[x_3, x_4, x_5, x_6, x_7]$.

\item[114.] $\mathfrak{g}_{7,1.7} \cong (247O)$\\
$[x_1,x_2] = x_4$, $[x_1, x_3] = x_5$, $[x_1, x_4] = x_6$, $[x_1,x_5] = x_7$, $[x_2, x_3] = x_6$,\\
$[x_2, x_4] = x_7$.\\
$i = 3$, $r= 1$, $c=5$, $F = \langle x_3, x_4, x_5, x_6, x_7\rangle = CPI$,\\
$Y = k[x_6, x_7, 2x_5x_6^2 - x_5^2x_7 - 2x_4x_6x_7 + 2x_3x_7^2]$, $M = k[x_3, x_4, x_5, x_6, x_7]$.

\item[115.] $\mathfrak{g}_{7,1.9} \cong (2457K)\cong R_{59}$\\
$[x_1,x_2] = x_4$, $[x_1, x_3] = x_6$, $[x_1, x_4] = x_5$, $[x_1,x_5] = x_7$, $[x_2, x_3] = x_7$,\\
$[x_2, x_4] = x_6$.\\
$i = 3$, $r= 1$, $c=5$, $F = \langle x_3, x_4, x_5, x_6, x_7\rangle = CPI$,\\
$Y = k[x_6, x_7, 2x_5x_6^2 - x_5^2x_7 -2x_3x_6x_7 + 2x_4x_7^2]$, $M = k[x_3, x_4, x_5, x_6, x_7]$.

\item[116.] $\mathfrak{g}_{7,2.36} \cong (257G)\cong R_{118}$\\
$[x_1,x_3] = x_5$, $[x_1, x_4] = x_7$, $[x_2, x_3] = x_6$, $[x_2,x_4] = -x_5$, $[x_3, x_6] = -x_7$.\\
$i = 3$, $r= 2$, $c=5$, $F = \langle x_1, x_2, x_5, x_6, x_7\rangle = CPI$,\\
$Y = k[x_5, x_7, 2x_5^2x_6 - 2x_1x_5x_7 + x_6^2x_7 - 2x_2x_7^2]$, $M = k[x_1, x_2, x_5, x_6, x_7]$.

\item[117.] $\mathfrak{g}_{7,2.1 (i_\lambda), \lambda = 0} \cong (2357B)\cong R_{75}$\\
$[x_1,x_2] = x_4$, $[x_1, x_3] = x_5$, $[x_1, x_4] = x_6$, $[x_1,x_6] = x_7$, $[x_2, x_3] = x_6$,\\
$[x_3,x_4] =-x_7$.\\
$i = 3$, $r= 2$, $c=5$, $F = \langle x_2, x_4, x_5, x_6, x_7\rangle = CPI$,\\
$Y = k[x_5, x_7, x_6^3 - 3x_4x_6x_7 + 3x_2x_7^2]$, $M = k[x_2, x_4, x_5, x_6, x_7]$.

\item[118.] $\mathfrak{g}_{7,1.13} \cong (23457E)\cong R_{30}$\\
$[x_1,x_2] = x_3$, $[x_1, x_3] = x_4$, $[x_1, x_4] = x_6$, $[x_1,x_5] = x_7$, $[x_1, x_6] = x_7$,\\
$[x_2, x_3] = x_5$, $[x_2, x_4] = x_7$.\\
$i = 3$, $r= 1$, $c=5$, $F = \langle x_3, x_4, x_5, x_6, x_7\rangle = CPI$,\\
$Y = k[x_5 -x_6, x_7, x_6^3 - 3x_5x_6^2 + 6x_4x_5x_7 - 6x_3x_7^2]$, $M = k[x_3, x_4, x_5, x_6, x_7]$.

\item[119.] $\mathfrak{g}_{7,2.24} \cong (2357A)\cong R_{73}$\\
$[x_1,x_2] = x_4$, $[x_1, x_4] = x_6$, $[x_1, x_5] = -x_7$, $[x_1,x_6] = x_7$, $[x_2, x_3] = x_5$,\\
$[x_3,x_4] = x_7$.\\
$i = 3$, $r= 2$, $c=5$, $F = \langle x_2, x_4, x_5, x_6, x_7\rangle = CPI$,\\
$Y = k[x_5 + x_6, x_7, x_6^3 + 3x_5x_6^2 - 6x_4x_5x_7 - 6x_2x_7^2]$, $M = k[x_2, x_4, x_5, x_6, x_7]$.

\item[120.] $\mathfrak{g}_{7,1.18} \cong (2457J)\cong R_{57}$\\
$[x_1,x_2] = x_4$, $[x_1, x_4] = x_5$, $[x_1, x_5] = x_7$, $[x_2,x_3] = x_6+x_7$, $[x_2, x_4] = x_6$.\\
$i = 3$, $r= 1$, $c=5$, $F = \langle x_3, x_4, x_5, x_6, x_7\rangle = CPI$,\\
$Y = k[x_6, x_7, x_5^2x_6 + x_5^2x_7 + 2x_3x_6x_7 - 2x_4x_6x_7 - 2x_4x_7^2]$, $M = k[x_3, x_4, x_5, x_6, x_7]$.

\item[121.] $\mathfrak{g}_{7,2.5} \cong (12457H)\cong R_{24}$\\
$[x_1,x_2] = x_3$, $[x_1, x_3] = x_4$, $[x_1, x_5] = x_6$, $[x_1,x_6] = x_7$, $[x_2, x_3] = x_5$,\\
$[x_2,x_4] = x_6$, $[x_3, x_4] = x_7$.\\
$i = 3$, $r= 2$, $c=5$, $F = \langle x_2, x_3, x_4, x_5, x_6, x_7\rangle$, no $CP$'s,\\
$Y = k[x_7, x_6^2 - 2x_5x_7, x_4x_5 - x_3x_6 + x_2x_7]$, $M = k[x_4, x_5, x_6, x_7, x_3x_6 - x_2x_7]$.

\item[122.] $\mathfrak{g}_{7,2.13} \cong (12457C)\cong R_{37}$\\
$[x_1,x_2] = x_4$, $[x_1, x_4] = x_5$, $[x_1, x_5] = x_6$, $[x_2,x_3] = x_6$, $[x_2, x_6] = x_7$,\\
$[x_4,x_5] = -x_7$.\\
$i = 3$, $r= 2$, $c=5$, $F = \langle x_1, x_3, x_4, x_5, x_6, x_7\rangle$, no $CP$'s,\\
$Y = k[x_7, x_6^2 - 2x_3x_7, x_5^2 - 2x_4x_6 - 2x_1x_7]$, $M = k[x_3, x_5, x_6, x_7, x_4x_6 + x_1x_7]$.

\item[123.] $\mathfrak{g}_{7,1.17} \cong (12457L)\cong R_{17}$\\
$[x_1,x_2] = x_3$, $[x_1, x_3] = x_4$, $[x_1, x_4] = x_6$, $[x_1,x_6] = x_7$, $[x_2, x_3] = x_5$, \\
$[x_2,x_5] = x_6$, $[x_2, x_6] = x_7$, $[x_3,x_4] = -x_7$, $[x_3, x_5] = x_7$.\\
$i = 3$, $r= 1$, $c=5$, $F = \langle x_1 - x_2, x_3, x_4, x_5, x_6, x_7\rangle$, no $CP$'s,\\
$Y = k[x_7, x_4^2 - 2x_3x_6 + x_5^2 - 2x_1x_7 + 2x_2x_7, x_6^2 - 2x_4x_7 - 2x_5x_7]$,\\
$M = k[x_4, x_5, x_6, x_7, x_3x_6 + x_1x_7 - x_2x_7]$.

\item[124.] $\mathfrak{g}_{7,2.1(v)} \cong (247Q)$\\
$[x_1,x_2] = x_4$, $[x_1, x_3] = x_5$, $[x_1, x_4] = x_6$, $[x_2,x_3] = x_6$, $[x_2, x_5] = x_7$,\\
$[x_3,x_4] = x_7$.\\
$i = 3$, $r= 2$, $c=5$, $F = \langle x_1, x_3, x_4, x_5, x_6, x_7\rangle$, no $CP$'s,\\
$Y = k[x_6, x_7, x_5x_6^2 - x_3x_6x_7 + x_4x_5x_7 + x_1x_7^2]$, $M = k[x_4, x_5, x_6, x_7, x_3x_6 - x_1x_7]$.

\item[125.] $\mathfrak{g}_{7,2.6} \cong (23457B) \cong R_{28}$\\
$[x_1,x_2] = x_3$, $[x_1, x_3] = x_4$, $[x_1, x_4] = x_6$, $[x_2,x_3] = x_5$, $[x_2, x_6] = x_7$,\\
$[x_3,x_4] = -x_7$.\\
$i = 3$, $r= 2$, $c=5$, $F = \langle x_1, x_3, x_4, x_5, x_6, x_7\rangle$, no $CP$'s,\\
$Y = k[x_5, x_7, x_5x_6^2 + x_4^2x_7 - 2x_3x_6x_7 - 2x_1x_7^2]$, $M = k[x_4, x_5, x_6, x_7, x_3x_6 + x_1x_7]$.

\item[126.] $\mathfrak{g}_{7,2.26} \cong (247J) \cong R_{85}$\\
$[x_1,x_2] = x_4$, $[x_1, x_3] = x_5$, $[x_1, x_5] = x_6$, $[x_2,x_5] = x_7$, $[x_3, x_4] = x_7$,\\
$[x_3,x_5] = x_6$.\\
$i = 3$, $r= 2$, $c=5$, $F = \langle x_1, x_2, x_4, x_5, x_6, x_7\rangle$, no $CP$'s,\\
$Y = k[x_6, x_7, x_4^2x_6 - 2x_4x_5x_7 + 2x_2x_6x_7 - 2x_1x_7^2]$,\\
$M = k[x_4, x_5, x_6, x_7, x_2x_6 - x_1x_7]$.

\item[127.] $\mathfrak{g}_{7,2.29} \cong (257L) \cong R_{104}$\\
$[x_1,x_2] = x_5$, $[x_1, x_5] = x_6$, $[x_2, x_3] = x_6$, $[x_2,x_5] = x_7$, $[x_3, x_4] = x_7$.\\
$i = 3$, $r= 2$, $c=5$, $F = \langle x_1, x_2, x_4, x_5, x_6, x_7\rangle$, no $CP$'s,\\
$Y = k[x_6, x_7, 2x_4x_6^2 - x_5^2x_7 + 2x_2x_6x_7 - 2x_1x_7^2]$, $M = k[x_4, x_5, x_6, x_7, x_2x_6 - x_1x_7]$.

\item[128.] $\mathfrak{g}_{7,2.34} \cong (247G) \cong R_{82}$\\
$[x_1,x_2] = x_4$, $[x_1, x_3] = x_5$, $[x_1, x_4] = x_7$, $[x_2,x_4] = x_6$, $[x_3, x_5] = x_7$.\\
$i = 3$, $r= 2$, $c=5$, $F = \langle x_1, x_2, x_4, x_5, x_6, x_7\rangle$, no $CP$'s,\\
$Y = k[x_6, x_7, x_5^2x_6 + x_4^2x_7 + 2x_1x_6x_7 - 2x_2x_7^2]$, $M = k[x_4, x_5, x_6, x_7, x_1x_6 - x_2x_7]$.

\item[129.] $\mathfrak{g}_{7,1.14} \cong (23457F) \cong R_{27}$\\
$[x_1,x_2] = x_3$, $[x_1, x_3] = x_4$, $[x_1, x_4] = x_6$, $[x_2,x_3] = x_5$, $[x_2, x_5] = -x_7$,\\
$[x_2,x_6] = -x_7$, $[x_3, x_4] = x_7$.\\
$i = 3$, $r= 1$, $c=5$, $F = \langle x_1, x_3, x_4, x_5, x_6, x_7\rangle$, no $CP$'s,\\
$Y = k[x_5- x_6, x_7, x_5^3 - 3x_5^2x_6 + 3x_4^2x_7 - 6x_3x_6x_7 + 6x_1x_7^2]$,\\
$M = k[x_4, x_5, x_6, x_7, x_3x_6 - x_1x_7]$.

\item[130.] $\mathfrak{g}_{7,2.17} \cong (2357C) \cong R_{56}$\\
$[x_1,x_2] = x_4$, $[x_1, x_4] = x_5$, $[x_1, x_5] = x_7$, $[x_2,x_3] = x_5$, $[x_2, x_4] = x_6$,\\
$[x_3,x_4] = -x_7$,.\
$i = 3$, $r= 2$, $c=5$, $F = \langle x_2, x_3, x_4, x_5, x_6, x_7\rangle$, no $CP$'s,\\
$Y = k[x_6, x_7, x_5^3 - 3x_4x_5x_7 + 3x_3x_6x_7 + 3x_2x_7^2]$, $M = k[x_4, x_5, x_6, x_7, x_3x_6 + x_2x_7]$.

\item[131.] $\mathfrak{g}_{7,1.2(iii)} \cong (2357D) \cong R_{54}$\\
$[x_1,x_2] = x_4$, $[x_1, x_3] = x_6$, $[x_1, x_4] = x_5$, $[x_1,x_5] = x_7$, $[x_2, x_3] = x_5$,\\
$[x_2,x_4] = x_6$, $[x_3, x_4] = -x_7$.\\
$i = 3$, $r= 1$, $c=5$, $F = \langle x_2, x_3, x_4, x_5, x_6, x_7\rangle$, no $CP$'s,\\
$Y = k[x_6, x_7, x_5^3 - 3x_5x_6^2 - 3x_4x_5x_7 + 3x_3x_6x_7 + 3x_2x_7^2]$,\\
$M = k[x_4, x_5, x_6, x_7, x_3x_6 + x_2x_7]$.

\item[132.] $\mathfrak{g}_{7,1.3(iv)} \cong (247R) \cong R_{87}$\\
$[x_1,x_2] = x_4$, $[x_1, x_3] = x_5$, $[x_1, x_4] = x_6$, $[x_2,x_3] = x_6$, $[x_2, x_4] = x_7$, \\
$[x_3,x_5] = x_7$.\\
$i = 3$, $r= 1$, $c=5$, $F = \langle x_1, x_2, x_4, x_5, x_6, x_7\rangle$, no $CP$'s,\\
$Y = k[x_6, x_7, 2x_5x_6^2 - x_4^2x_7 - x_5^2x_7 + 2x_2x_6x_7 - 2x_1x_7^2]$,\\
$M = k[x_4, x_5, x_6, x_7, x_2x_6 - x_1x_7]$.

\item[133.] $\mathfrak{g}_{7,1.5} \cong (23457D) \cong R_{12}$ (quasi quadratic)\\
$[x_1,x_2] = x_3$, $[x_1, x_3] = x_4$, $[x_1, x_4] = x_5$, $[x_1,x_5] = x_6$, $[x_2, x_3] = x_6$,\\
$[x_2,x_5] = -x_7$, $[x_3, x_4] = x_7$.\\
$i = 3$, $r= 1$, $c=5$, $F = \mathfrak{g}_{7,1.5}$, no $CP$'s,\\
$Y = k[x_6, x_7, 2x_4x_6^2 - x_5^2x_6 + x_4^2x_7- 2x_3x_5x_7 + 2x_2x_6x_7 + 2x_1x_7^2]$,\\
$M = k[x_4, x_5, x_6, x_7, -x_3x_5 + x_2x_6 + x_1x_7]$.

\item[134.] $\mathfrak{g}_{7,2.35} \cong (247K) \cong R_{84}$ (quasi quadratic)\\
$[x_1,x_2] = x_4$, $[x_1, x_3] = x_5$, $[x_1, x_5] = x_6$, $[x_2,x_4] = x_6$, $[x_2, x_5] = x_7$,\\
$[x_3,x_4] = x_7$.\\
$i = 3$, $r= 2$, $c=5$, $F = \mathfrak{g}_{7,2.35}$, no $CP$'s,\\
$Y = k[x_6, x_7, 2x_3x_6^2 - x_5^2x_6 + 2x_4x_5x_7 - 2x_2x_6x_7 + 2x_1x_7^2]$,\\
$M = k[x_4, x_5, x_6, x_7, x_3x_6^2 - x_2x_6x_7 + x_1x_7^2]$,\\
$M_1 = k[x_6, x_7, x_5^2 - 2x_3x_6 + x_2x_7, 2x_4x_5 - x_2x_6 + 2x_1x_7, x_5x_6 + x_4x_7]$.

\item[135.] $\mathfrak{g}_{7,1.19} \cong (247H) \cong R_{81}$ (quasi quadratic)\\
$[x_1,x_2] = x_4$, $[x_1, x_3] = x_5$, $[x_1, x_4] = x_7$, $[x_1,x_5] = x_6$, $[x_2, x_4] = x_6$,\\
$[x_3,x_5] = x_7$.\\
$i = 3$, $r= 1$, $c=5$, $F = \mathfrak{g}_{7,1.19}$, no $CP$'s,\\
$Y = k[x_6, x_7, 2x_3x_6^2 - x_5^2x_6 - x_4^2x_7- 2x_1x_6x_7 + 2x_2x_7^2]$,\\
$M = k[x_4, x_5, x_6, x_7, x_3x_6^2 - x_1x_6x_7 + x_2x_7^2]$.\\
$M_1 = k[x_6, x_7, x_5^2 - 2x_3x_6 + x_1x_7, x_4^2 + x_1x_6 - 2x_2x_7, x_4x_6 + x_5x_7]$.

\item[136.] $\mathfrak{g}_{7,1.1(iii)} \cong (23457G) \cong R_{11}$ (quasi quadratic, see Example 33)\\
$[x_1,x_2] = x_3$, $[x_1, x_3] = x_4$, $[x_1, x_4] = x_5$, $[x_1,x_5] = x_6$, $[x_2, x_3] = x_5$,\\
$[x_2,x_4] = x_6, [x_2, x_5] = -x_7, [x_3, x_4] = x_7$.\\
$i = 3$, $r= 1$, $c=5$, $cod = 3$, $F = \mathfrak{g}_{7,1.1(iii)}$, no $CP$'s,\\
$Y = k[x_6, x_7, 2x_5^3 - 3x_4^2x_7 - 6x_4x_5x_6 + 6x_3x_6^2 + 6x_3x_5x_7 - 6x_2x_6x_7 - 6x_1x_7^2]$,\\
$M = k[x_4, x_5, x_6, x_7, x_3x_6^2 + x_3x_5x_7 - x_2x_6x_7 - x_1x_7^2]$.\\
$Y_\xi = k [x_5, x_6, x_7, x_4x_6 - x_3x_7, x_4^2x_7 - 2x_3x_6^2 + 2x_2x_6x_7 + 2x_1x_7^2]$, where $\xi = x_5^\ast$,\\
$M_1 = k[x_5, x_6, x_7, x_3x_6 - x_4x_5, x_4^2 + 2x_1x_7 + 2x_2x_6 - 2x_3x_5]$.
\end{itemize}

{\bf II. $\mathfrak{g}$ is not coregular} (23 cases) \\
\begin{itemize}
\item[137.] $\mathfrak{g}_{7,1.6} \cong (123457B) \cong R_{9}$\\
$[x_1,x_2] = x_3$, $[x_1, x_3] = x_4$, $[x_1, x_4] = x_5$, $[x_1,x_5] = x_6$, $[x_1, x_6] = x_7$,\\
$[x_2,x_3] = x_7$.\\
$i = 3$, $r= 1$, $c=5$, $F = \langle x_4, x_5, x_6, x_7\rangle$, $CPI = \langle x_3, x_4, x_5, x_6, x_7\rangle$,\\
$cod = 1$, $p = x_7$, $Y = k[x_7, f, g, h]$, $f = x_6^2 - 2x_5x_7$, $g = x_6^3 - 3x_5x_6x_7 + 3x_4x_7^2$, $h = (f^3-g^2) / x_7^2$, $Q(Y) = k(x_7, f, g)$,\\
relation: $f^3 - g^2 - x_7^2h = 0$, $M = k[x_3, x_4, x_5, x_6, x_7]$.

\item[138.] $\mathfrak{g}_{7,1.15} \cong (13457B) \cong R_{50}$\\
$[x_1,x_2] = x_4$, $[x_1, x_4] = x_5$, $[x_1, x_5] = x_6$, $[x_1,x_6] = x_7$, $[x_2, x_3] = x_7$,\\
$[x_2,x_4] = x_7$.\\
$i = 3$, $r= 1$, $c=5$, $F = \langle x_3-x_4, x_5, x_6, x_7\rangle$, $CPI = \langle x_3, x_4, x_5, x_6, x_7\rangle$,\\
$cod = 1$, $p = x_7$, $Y = k[x_7, f, g, h]$, $f = x_6^2 - 2x_5x_7$,\\
$g = x_6^3 - 3x_5x_6x_7 + 3x_4x_7^2- 3x_3x_7^2$, $h = (f^3-g^2) / x_7^2$, \\
relation: $f^3 - g^2 - x_7^2h = 0$, $Q(Y) = k(x_7, f, g)$, $M = k[x_3, x_4, x_5, x_6, x_7]$.

\item[139.] $\mathfrak{g}_{7,2.16} \cong (13457A) \cong R_{51}$\\
$[x_1,x_2] = x_4$, $[x_1, x_4] = x_5$, $[x_1, x_5] = x_6$, $[x_1,x_6] = x_7$, $[x_2, x_3] = x_7$.\\
$i = 3$, $r= 2$, $c=5$, $F = \langle x_4, x_5, x_6, x_7\rangle$, $CPI = \langle x_3, x_4, x_5, x_6, x_7\rangle$,\\
$cod = 1$, $p = x_7$, $Y = k[x_7, f, g, h]$, $f = x_6^2 - 2x_5x_7$, $g = x_6^3 - 3x_5x_6x_7 + 3x_4x_7^2$, $h = (f^3-g^2) / x_7^2$, \\
relation: $f^3 - g^2 - x_7^2h = 0$, $Q(Y) = k(x_7, f, g)$, $M = k[x_3, x_4, x_5, x_6, x_7]$.

\item[140.] $\mathfrak{g}_{7,2.32} \cong (1357E) \cong R_{67}$\\
$[x_1,x_2] = x_4$, $[x_1, x_3] = x_5$, $[x_1, x_4] = x_6$, $[x_1,x_6] = x_7$, $[x_3, x_5] = x_7$.\\
$i = 3$, $r= 2$, $c=5$, $F = \langle x_2, x_4, x_6, x_7\rangle$, $CPI = \langle x_2, x_4, x_5, x_6, x_7\rangle$,\\
$cod = 1$, $p = x_7$, $Y = k[x_7, f, g, h]$, $f = x_6^2 - 2x_4x_7$, $g = x_6^3 - 3x_4x_6x_7 + 3x_2x_7^2$, $h = (f^3-g^2) / x_7^2$, \\
relation: $f^3 - g^2 - x_7^2h = 0$, $Q(Y) = k(x_7, f, g)$, $M = k[x_2, x_4, x_5, x_6, x_7]$.

\item[141.] $\mathfrak{g}_{7,0.3} \cong (123457E) \cong R_{7}$\\
$[x_1,x_2] = x_3$, $[x_1, x_3] = x_4$, $[x_1, x_4] = x_5$, $[x_1,x_5] = x_6$, $[x_1, x_6] = x_7,$\\
$[x_2, x_3] = x_6+x_7, [x_2, x_4] = x_7$.\\
$i = 3$, $r= 0$, $c=5$, $cod = 2$, $F = \langle x_3, x_4, x_5, x_6, x_7\rangle = CPI$, \\
$Y = k[x_7, f, g, h]$, $f = x_6^2 - 2x_5x_7$,\\
$g = 2x_6^3 - 3x_5^2x_7 + 6x_4x_6x_7 -6x_5x_6x_7 - 6x_3x_7^2 + 6x_4x_7^2$, $h = (4f^3-g^2) / x_7$, \\
relation: $4f^3 - g^2 - x_7h = 0$, $Q(Y) = k(x_7, f, g)$, $M = k[x_3, x_4, x_5, x_6, x_7]$.

\item[142.] $\mathfrak{g}_{7,1.01(ii)} \cong (12457B) \cong R_{47}$\\
$[x_1,x_2] = x_4$, $[x_1, x_4] = x_5$, $[x_1, x_5] = x_6$, $[x_1,x_6] = x_7$,\\
$[x_2, x_3] = x_6 + x_7, [x_3, x_4] = -x_7$.\\
$i = 3$, $r= 1$, $c=5$, $cod = 2$, $F = \langle x_2, x_4, x_5, x_6, x_7\rangle = CPI$, \\
$Y = k[x_7, f, g, h]$, $f = x_6^2 - 2x_5x_7$,\\
$g = 2x_6^3 - 3x_5^2x_7 + 6x_4x_6x_7 -6x_5x_6x_7 - 6x_2x_7^2 + 6x_4x_7^2$, $h = (4f^3-g^2) / x_7$, \\
relation: $4f^3 - g^2 - x_7h = 0$, $Q(Y) = k(x_7, f, g)$, $M = k[x_2, x_4, x_5, x_6, x_7]$.

\item[143.] $\mathfrak{g}_{7,1.3(i_\lambda),\lambda=0} \cong (1357N) (\xi = 0) \cong R_{62}^0$\\
$[x_1,x_2] = x_4$, $[x_1, x_3] = x_5$, $[x_1, x_4] = x_6$, $[x_1,x_6] = x_7$, $[x_2, x_3] = x_6,$\\
$[x_2, x_5] = x_7, [x_3, x_5] = x_7$.\\
$i = 3$, $r= 1$, $c=5$, $cod = 2$, $F = \langle x_2-x_3, x_4, x_5, x_6, x_7\rangle = CPI$, \\
$Y = k[x_7, f, g, h]$, $f = x_6^2 - 2x_4x_7$, $g = x_6^3 - 3x_4x_6x_7 + 3x_5x_6x_7 +3x_2x_7^2 - 3x_3x_7^2$, $h = (f^3-g^2) / x_7$, \\
relation: $f^3 - g^2 - x_7h = 0$, $Q(Y) = k(x_7, f, g)$, $M = k[x_2-x_3, x_4, x_5, x_6, x_7]$.

\item[144.] $\mathfrak{g}_{7,1.10} \cong (13457F) \cong R_{29}$\\
$[x_1,x_2] = x_3$, $[x_1, x_3] = x_4$, $[x_1, x_4] = x_6$, $[x_1,x_6] = x_7$, $[x_2, x_3] = x_5,$ \\
$[x_2, x_5] = x_7$.\\
$i = 3$, $r= 1$, $c=5$, $cod = 2$, $F = \langle x_3, x_4, x_5, x_6, x_7\rangle = CPI$, \\
$Y = k[x_7, f, g, h]$, $f = x_6^2 - 2x_4x_7$,
$g = 2x_6^3 - 3x_5^2x_7 - 6x_4x_6x_7 +6x_3x_7^2$,\\
$h = (4f^3-g^2) / x_7$, \\
relation: $4f^3 - g^2 - x_7h = 0$, $Q(Y) = k(x_7, f, g)$, $M = k[x_3, x_4, x_5, x_6, x_7]$.

\item[145.] $\mathfrak{g}_{7,2.33} \cong (1357I) \cong R_{71}$\\
$[x_1,x_2] = x_4$, $[x_1, x_4] = x_6$, $[x_1, x_6] = x_7$, $[x_2,x_3] = x_5$, $[x_3, x_5] = x_7$.\\
$i = 3$, $r= 2$, $c=5$, $cod = 2$, $F = \langle x_2, x_4, x_5, x_6, x_7\rangle = CPI$, \\
$Y = k[x_7, f, g, h]$, $f = x_6^2 - 2x_4x_7$, $g = 2x_6^3 + 3x_5^2x_7 - 6x_4x_6x_7 + 6x_2x_7^2$,\\
$h = (4f^3-g^2) / x_7$, \\
relation: $4f^3 - g^2 - x_7h = 0$, $Q(Y) = k(x_7, f, g)$, $M = k[x_2, x_4, x_5, x_6, x_7]$.

\item[146.] $\mathfrak{g}_{7,1.12} \cong (13457D) \cong R_{49}$\\
$[x_1,x_2] = x_4$, $[x_1, x_4] = x_5$, $[x_1, x_5] = x_6$, $[x_1,x_6] = x_7$, $[x_2, x_3] = x_7,$ \\
$[x_2, x_4] = x_6, [x_2, x_5] = x_7$.\\
$i = 3$, $r= 1$, $c=5$, $cod = 2$, $F = \langle x_3, x_4, x_5, x_6, x_7\rangle = CPI$, \\
$Y = k[x_7, f, g, h]$, $f = x_6^2 + 2x_3x_7 - 2x_5x_7$, $g = x_6^3 - 3x_5x_6x_7 +3x_4x_7^2$,\\
$h = (f^3-g^2) / x_7$, \\
relation: $f^3 - g^2 - x_7h = 0$, $Q(Y) = k(x_7, f, g)$, $M = k[x_3, x_4, x_5, x_6, x_7]$.

\item[147.] $\mathfrak{g}_{7,2.25} \cong (1357B) \cong R_{100}$\\
$[x_1,x_2] = x_5$, $[x_1, x_5] = x_6$, $[x_1, x_6] = x_7$, $[x_2,x_3] = x_6$, $[x_3, x_4] = -x_7,$\\
$[x_3, x_5] = -x_7$.\\
$i = 3$, $r= 2$, $c=5$, $cod = 2$, $F = \langle x_2, x_4, x_5, x_6, x_7\rangle = CPI$, \\
$Y = k[x_7, f, g, h]$, $f = x_6^2 + 2x_4x_7 - 2x_5x_7$, $g = x_6^3 - 3x_5x_6x_7 + 3x_2x_7^2$,\\
$h = (f^3-g^2) / x_7$, \\
relation: $f^3 - g^2 - x_7h = 0$, $Q(Y) = k(x_7, f, g)$, $M = k[x_2, x_4, x_5, x_6, x_7]$.

\item[148.] $\mathfrak{g}_{7,0.5} \cong (12457N)(\xi = 1) \cong R_{16}$\\
$[x_1,x_2] = x_3$, $[x_1, x_3] = x_4$, $[x_1, x_4] = x_6 + x_7$, $[x_1,x_6] = x_7$, \\
$[x_2, x_3] = x_5, [x_2, x_5] = x_6, [x_3, x_5] = x_7$.\\
$i = 3$, $r= 0$, $c=5$, $cod = 2$, $F = \langle x_2, x_3, x_4, x_5, x_6, x_7\rangle$, no $CP$'s, \\
$Y = k[x_7, f, g, h]$, $f = x_6^2 - 2x_4x_7 + 2x_6x_7$,\\
$g = 2x_6^3 + 3x_6^2x_7 + 3x_4^2x_7 +3x_5^2x_7 - 6x_3x_6x_7 - 6x_4x_6x_7 + 6x_2x_7^2$,\\
$h = (4f^3-g^2) / x_7$, \\
relation: $4f^3 - g^2 - x_7h = 0$, $Q(Y) = k(x_7, f, g)$,\\
$M = k[x_4, x_5, x_6, x_7, x_3x_6 - x_2x_7]$.

\item[149.] $\mathfrak{g}_{7,1.02} \cong (12457K) \cong R_{23}$\\
$[x_1,x_2] = x_3$, $[x_1, x_3] = x_4$, $[x_1, x_5] = x_6$, $[x_2,x_3] = x_5$, $[x_2, x_4] = x_6,$\\
$[x_2, x_5] = x_7, [x_2, x_6] = x_7, [x_3,x_5] = -x_7$.\\
$i = 3$, $r= 1$, $c=5$, $cod = 2$, $F = \langle x_1, x_3, x_4, x_5, x_6, x_7\rangle$, no $CP$'s, \\
$Y = k[x_7, f, g, h]$, $f = x_6^2 - 2x_4x_7$, $g = x_6^3 + 3x_4x_5x_7 - 3x_3x_6x_7 -3x_4x_6x_7 - 3x_1x_7^2$, $h = (f^3-g^2) / x_7$, \\
relation: $f^3 - g^2 - x_7h = 0$, $Q(Y) = k(x_7, f, g)$,\\
$M = k[x_4, x_5, x_6, x_7, x_3x_6 + x_1x_7]$.

\item[150.] $\mathfrak{g}_{7,1.21} \cong (12457F) \cong R_{36}$ (see Example 35)\\
$[x_1,x_2] = x_4$, $[x_1, x_4] = x_5$, $[x_1, x_5] = x_6$, $[x_2,x_3] = x_6$, $[x_2, x_4] = x_6,$ \\
$[x_2, x_6] = x_7, [x_4, x_5] = -x_7$.\\
$i = 3$, $r= 1$, $c=5$, $cod = 2$, $F = \langle x_1, x_3, x_4, x_5, x_6, x_7\rangle$, no $CP$'s, \\
$Y = k[x_7, f, g, h]$, $f = x_6^2 - 2x_3x_7$, $g = 2x_6^3 + 3x_5^2x_7 - 6x_4x_6x_7 - 6x_1x_7^2$,\\
$h = (4f^3-g^2) / x_7$, \\
relation: $4f^3 - g^2 - x_7h = 0$, $Q(Y) = k(x_7, f, g)$,\\
$M = k[x_3, x_5, x_6, x_7, x_4x_6 + x_1x_7]$.

\item[151.] $\mathfrak{g}_{7,1.1(i_\lambda),\lambda =1} \cong (123457G) \cong R_{1}^1$\\
$[x_1,x_2] = x_3$, $[x_1, x_3] = x_4$, $[x_1, x_4] = x_5$, $[x_1,x_5] = x_6$, $[x_1, x_6] = x_7,$\\
$[x_2, x_3] = x_5, [x_2, x_4] = x_6, [x_2,x_5] = x_7$.\\
$i = 3$, $r= 1$, $c=5$, $cod = 3$, $F = \langle x_3, x_4, x_5, x_6, x_7\rangle = CPI$, \\
$Y = k[x_7, f, g, h]$, $f = x_6^3 - 3x_5x_6x_7 + 3x_4x_7^2$,\\
$g = x_6^4 - 4x_5x_6^2x_7 + 2x_5^2x_7^2 + 4x_4x_6x_7^2 - 4x_3x_7^3$, $h = (f^4-g^3) / x_7^3$, \\
relation: $f^4 - g^3 - x_7^3h = 0$, $Q(Y) = k(x_7, f, g)$, $M = k[x_3, x_4, x_5, x_6, x_7]$.

\item[152.] $\mathfrak{g}_{7,2.14} \cong (12357A) \cong R_{45}$\\
$[x_1,x_2] = x_4$, $[x_1, x_4] = x_5$, $[x_1, x_5] = x_6$, $[x_1,x_6] = x_7$, $[x_2, x_3] = x_5,$ \\
$[x_3, x_4] = -x_6, [x_3, x_5] = -x_7$.\\
$i = 3$, $r= 2$, $c=5$, $cod = 3$, $F = \langle x_2, x_4, x_5, x_6, x_7\rangle = CPI$, \\
$Y = k[x_7, f, g, h]$, $f = x_6^3 - 3x_5x_6x_7 + 3x_4x_7^2$, \\
$g = x_6^4 - 4x_5x_6^2x_7 + 4x_4x_6x_7^2 + 2x_5^2x_7^2 - 4x_2x_7^3$, $h = (f^4-g^3) / x_7^3$, \\
relation: $f^4 - g^3 - x_7^3h = 0$, $Q(Y) = k(x_7, f, g)$, $M = k[x_2, x_4, x_5, x_6, x_7]$.

\item[153.] $\mathfrak{g}_{7,0.2} \cong (123457H) \cong R_{4}$\\
$[x_1,x_2] = x_3$, $[x_1, x_3] = x_4$, $[x_1, x_4] = x_5$, $[x_1,x_5] = x_6$, $[x_1, x_6] = x_7,$ \\
$[x_2, x_3] = x_5+x_7, [x_2, x_4] = x_6, [x_2, x_5] = x_7$.\\
$i = 3$, $r= 0$, $c=5$, $cod = 3$, $F = \langle x_3, x_4, x_5, x_6, x_7\rangle = CPI$,\\
$Y = k[x_7, f, g, h]$, $f = x_6^3 - 3x_5x_6x_7 + 3x_4x_7^2$,\\
$g = x_6^4 - 4x_5x_6^2x_7 + 2x_5^2x_7^2 + 4x_4x_6x_7^2 - 2x_6^2x_7^2 - 4x_3x_7^3 + 4x_5x_7^3$,\\
$h = (f^4-g^3-6x_7^2f^2g) / x_7^3$, relation: $f^4 - g^3 - 6x_7^2f^2g - x_7^3 h = 0$,\\
$Q(Y) = k(x_7, f, g)$, $M = k[x_3, x_4, x_5, x_6, x_7]$.

\item[154.] $\mathfrak{g}_{7,1.01(i)} \cong (12357B) \cong R_{43}$ (see Example 34)\\
$[x_1,x_2] = x_4$, $[x_1, x_4] = x_5$, $[x_1, x_5] = x_6$, $[x_1,x_6] = x_7$, $[x_2, x_3] = x_5 + x_7,$\\
$[x_3, x_4] = -x_6, [x_3, x_5] = -x_7$.\\
$i = 3$, $r= 1$, $c=5$, $cod = 3$, $F = \langle x_2, x_4, x_5, x_6, x_7\rangle = CPI$, \\
$Y = k[x_7, f, g, h]$, $f = x_6^3 - 3x_5x_6x_7 + 3x_4x_7^2$,\\
$g = x_6^4 - 4x_5x_6^2x_7 - 2x_6^2x_7^2 + 2x_5^2x_7^2 + 4x_4x_6x_7^2 + 4x_5x_7^3 - 4x_2x_7^3$,\\
$h = (f^4-g^3 - 6x_7^2f^2g) / x_7^3$, relation: $f^4 - g^3 - 6x_7^2f^2g - x_7^3 h = 0$, \\
$Q(Y) = k(x_7, f, g)$, $M = k[x_2, x_4, x_5, x_6, x_7]$.

\item[155.] $\mathfrak{g}_{7,1.1(i_\lambda), \lambda = 0} \cong (123457I)(\xi = 0) \cong R_{1}^0$\\
$[x_1,x_2] = x_3$, $[x_1, x_3] = x_4$, $[x_1, x_4] = x_5$, $[x_1,x_5] = x_6$, $[x_1, x_6] = x_7,$\\
$[x_2, x_3] = x_5, [x_2, x_4] = x_6, [x_3, x_4] = x_7$.\\
$i = 3$, $r= 1$, $c=5$, $cod = 3$, $F = \langle x_2, x_3, x_4, x_5, x_6, x_7\rangle$, no $CP$'s, \\
$Y = k[x_7, f, g, h]$, $f = x_6^2 - 2x_5x_7$,\\
$g = 2x_6^5 - 10x_5x_6^3 x_7 + 15x_5^2x_6x_7^2 - 15x_4x_5x_7^3 + 15x_3x_6x_7^3 - 15x_2x_7^4$,\\
$h = (4f^5-g^2) / x_7^3$, relation: $4f^5 - g^2 - x_7^3h = 0$, $Q(Y) = k(x_7, f, g)$, \\
$M = k[x_4, x_5, x_6, x_7, x_3x_6 - x_2x_7]$.

\item[156.] $\mathfrak{g}_{7,4.2} \cong (37A) \cong R_{127}$\\
$[x_1,x_2] = x_5$, $[x_1, x_3] = x_6$, $[x_1, x_4] = x_7$.\\
$i = 5$, $r= 4$, $c=6$, $cod = 3$, $F = \langle x_2, x_3, x_4, x_5, x_6, x_7\rangle = CPI$, \\
$Y = k[x_5, x_6, x_7, f, g, h]$, $f = x_3x_5 - x_2x_6$, $g = x_4x_6 - x_3x_7$, $h = x_4x_5 - x_2x_7$, \\
relation: $x_7f + x_5g - x_6h = 0$, $Q(Y) = k(x_5, x_6, x_7, f, g)$, $M = k[x_2, x_3, x_4, x_5, x_6, x_7]$.

\item[157.] $\mathfrak{g}_{7,3.20} \cong (247A)$\\
$[x_1,x_2] = x_4$, $[x_1, x_3] = x_5$, $[x_1, x_4] = x_6, [x_1, x_5] = x_7$.\\
$i = 5$, $r= 3$, $c=6$, $cod = 4$, $F = \langle x_2, x_3, x_4, x_5, x_6, x_7\rangle = CPI$, \\
$Y = k[x_6, x_7, f_1, f_2, f_3, f_4]$, $f_1 = x_4^2 - 2x_2x_6$, $f_2 = x_5^2 - 2x_3x_7$, $f_3 = x_5x_6 - x_4x_7$, $f_4 = -x_4x_5 + x_3x_6 + x_2x_7$,\\
relation: $x_7^2f_1 + x_6^2f_2 - f_3^2 + 2x_6x_7f_4 = 0$, $Q(Y) = k(x_6, x_7, f_1, f_2, f_3)$,\\
$M = k[x_2, x_3, x_4, x_5, x_6, x_7]$.

\item[158.] $\mathfrak{g}_{7,3.2} \cong (2457A) \cong R_{78}$\\
$[x_1,x_2] = x_4$, $[x_1, x_3] = x_5$, $[x_1, x_4] = x_6, [x_1, x_6] = x_7$.\\
$i = 5$, $r= 3$, $c=6$, $cod = 4$, $F = \langle x_2, x_3, x_4, x_5, x_6, x_7\rangle = CPI$, \\
$Y = k[f_1, f_2, \ldots, f_{13}]$, $Q(Y) = k(f_1, f_2, f_3, f_4, f_5)$, $M = k[x_2, x_3, x_4, x_5, x_6, x_7]$ where the following generators of $Y$ satisfy 63 relations (which will be omitted).\\
$f_1 = x_7$, $f_2 = x_5$, $f_3 = x_6^2 - 2x_4x_7$, $f_4 = x_5x_6 - x_3x_7$,\\
$f_5 = x_6^3 - 3x_4x_6x_7 + 3x_2x_7^2$, $f_6 = x_4x_5x_6 - x_3x_6^2 + 2x_3x_4x_7 - 3x_2x_5x_7$,\\
$f_7 = 2x_4x_5^2 - 2x_3x_5x_6 + x_3^2x_7$, \\
$f_8 = x_4x_5x_6^2 - x_3x_6^3 - 4x_4^2x_5x_7 + 3x_3x_4x_6x_7 + 3x_2x_5x_6x_7 - 3x_2x_3x_7^2$,\\
$f_9 = 3x_4^2x_6^2 - 6x_2x_6^3 - 8x_4^3x_7 + 18x_2x_4x_6x_7 - 9x_2^2x_7^2$,\\
$f_{10} = 4x_4^2 x_5^2 - 2x_3 x_4x_5x_6 - 6x_2x_5^2x_6 + x_3^2x_6^2 - 2x_3^2x_4x_7 + 6x_2x_3x_5x_7$,\\
$f_{11} = 6x_3x_4x_5^2 - 6x_2x_5^3 - 3x_3^2x_5x_6 + x_3^3x_7$,\\
$f_{12} = 2x_4^2x_5^2x_6 + 2x_3x_4x_5x_6^2 - 6x_2x_5^2x_6^2 - x_3^2x_6^3 - 8x_3x_4^2x_5x_7 + 6x_2x_4x_5^2x_7 + 3x_3^2x_4x_6x_7 + 6x_2 x_3x_5x_6x_7 - 3x_2x_3^2x_7^2$,\\
$f_{13} = 8x_4^3x_5^3 - 6x_3x_4^2x_5^2x_6 - 18x_2x_4x_5^3 x_6 - 3x_3^2 x_4x_5x_6^2 + 18x_2x_3x_5^2x_6^2 + x_3^3 x_6^3 + 12x_3^2 x_4^2x_5x_7 -18x_2x_3x_4x_5^2x_7 + 18x_2^2x_5^3x_7 - 3x_3^3x_4x_6x_7 - 9x_2x_3^2x_5x_6x_7 + 3x_2x_3^3x_7^2$.

\item[159.] $\mathfrak{g}_{7,2.3} \cong (123457A) \cong R_{10}$ (the 7-dim. standard filiform)\\
$[x_1,x_2] = x_3$, $[x_1, x_3] = x_4$, $[x_1, x_4] = x_5, [x_1, x_5] = x_6, [x_1, x_6] = x_7$.\\
$i = 5$, $r= 2$, $c=6$, $cod = 5$, $F = \langle x_2, x_3, x_4, x_5, x_6, x_7\rangle = CPI$, \\
$Y = k[f_1, f_2, \ldots, f_{23}]$ (see Example 28),\\
$M = k[x_2, x_3, x_4, x_5, x_6, x_7]$, $Q(Y) = k(f_1, f_2, f_3, f_4, f_5)$, where $f_1 = x_7$,\\
$f_2 = x_6^2 - 2x_5x_7$, $f_3 = x_6^3 - 3x_5x_6x_7 + 3x_4x_7^2$, $f_4 = x_5^2 - 2x_4x_6 + 2x_3x_7$,\\
$f_5 = 2x_4x_6^2 - x_5^2x_6 + x_4x_5x_7 - 5x_3x_6x_7 + 5x_2x_7^2$ are algebraically independent over $k$.
\end{itemize}

{\bf 3.3 A strongly complete subalgebra $M$ for the Lie algebras 1-82 of [O5,5]}\\
\ \\
For the majority of these Lie algebras a $CPI$, say $\mathfrak{h}$, was provided in [O5,5].  Then $M = S(\mathfrak{h})$ has all the required properties (see (1) of Examples 19).  Therefore we only have to consider the 6 remaining cases (i.e. without $CP$'s) listed below.\\
\ \\
{\underline{$\dim \mathfrak{g} = 5$}}
\begin{itemize}
\item[7.] $\mathfrak{g}_{5,4}$ (quadratic)\\
$[x_1, x_2] = x_3$, $[x_1, x_3] = x_4$, $[x_2, x_3] = x_5$.\\
$i=3$, $r=2$, $c=4$, $cod =3$, $F = \mathfrak{g}_{5,4}$,\\
$Y = k[x_4, x_5, x_3^2 + 2x_1x_5 - 2x_2x_4]$, $M = k[x_3, x_4, x_5, x_2x_4 - x_1x_5] = Y_\xi$, $\xi = x_3^\ast$.
\end{itemize}

{\underline{$\dim \mathfrak{g} = 6$}}
\begin{itemize}
\item[21.] $\mathfrak{g}_{6,18} \cong M_{21} \cong C_2$\\
$[x_1, x_2] = x_3$, $[x_1, x_3] = x_4$, $[x_1, x_4] = x_5$, $[x_2, x_5] = x_6$, $[x_3, x_4] = -x_6$.\\
$i=2$, $r=2$, $c=4$, $cod =1$, $p=x_6$, $F = \langle x_1, x_3, x_4, x_5, x_6\rangle$,\\
$Y = k[x_6, x_4^2 -2x_3x_5 - 2x_1x_6]$, $M = k[x_4, x_5, x_6, x_3x_5 + x_1x_6]$.

\item[22.] $\mathfrak{g}_{6,3} \cong M_{3} \cong C_{21}$ (quadratic)\\
$[x_1, x_2] = x_4$, $[x_1, x_3] = x_5$, $[x_2, x_3] = x_6$.\\
$i=4$, $r=3$, $c=5$, $cod =3$, $F = \mathfrak{g}_{6,3}$,\\
$Y = k[x_4, x_5, x_6, x_3x_4 - x_2x_5 + x_1x_6]$,\\
$M = k[x_3, x_4, x_5, x_6, x_2x_5 - x_1x_6] = Y_\xi$, $\xi = x_4^\ast$.

\item[28.] $\mathfrak{g}_{6,20} \cong M_{22} \cong C_1$\\
$[x_1, x_2] = x_3$, $[x_1, x_3] = x_4$, $[x_1, x_4] = x_5$, $[x_2, x_3] = x_5$, $[x_2, x_5] = x_6,$\\
$[x_3, x_4] = -x_6$.\\
$i=2$, $r=1$, $c=4$, $cod =2$, $F = \langle x_1, x_3, x_4, x_5, x_6\rangle$,\\
$Y = k[x_6, 2x_5^3 + 3x_4^2x_6 - 6x_3x_5x_6 - 6x_1x_6^2]$, $M = k[x_4, x_5, x_6, x_3x_5 + x_1x_6]$.
\end{itemize}

{\underline{$\dim \mathfrak{g} = 7$}}
\begin{itemize}
\item[69.] $\mathfrak{g}_{7, 3.13} \cong (257K) \cong R_{105}$\\
$[x_1, x_2] = x_5$, $[x_1, x_5] = x_6$, $[x_2, x_5] = x_7$, $[x_3, x_4] = x_7$.\\
$i=3$, $r=3$, $c=5$, $cod =1$, $p=x_7$, $F = \langle x_1, x_2, x_5, x_6, x_7\rangle$,\\
$Y = k[x_6, x_7, x_5^2 -2x_2x_6 + 2x_1x_7]$, $M = k[x_4, x_5, x_6, x_7, x_2x_6 - x_1x_7]$.

\item[77.] $\mathfrak{g}_{7, 3.22} \cong (247D)$\\
$[x_1, x_2] = x_4$, $[x_1, x_3] = x_5$, $[x_1, x_5] = x_7$, $[x_2, x_5] = x_6, [x_3, x_4] = x_6$.\\
$i=3$, $r=3$, $c=5$, $cod =1$, $p=x_6$, $F = \langle x_1, x_2, x_4, x_5, x_6, x_7\rangle$,\\
$Y = k[x_6, x_7, x_4x_5 + x_1x_6 - x_2x_7]$, $M = k[x_4, x_5, x_6, x_7, x_1x_6 - x_2x_7]$.
\end{itemize}

{\bf 3.4. A few observations}\\
\\
{\it 3.4.1. Lie algebras having isomorphic Poisson centers}\\
By considering $Y(\mathfrak{g})$ as the subalgebra of all elements of $k[x_1,\ldots, x_n]$ annihilated by $\sum\limits_{j=1}^n [x_i, x_j] \frac{\partial}{\partial x_j}$, $i : 1,\ldots,n$, sometimes one can reduce its determination to a former case.  For example: 83, 84, 85, 86, 87; 89, 90, 91, 92; 93, 94, 95, 96; 98, 99; 103, 104; 106, 107; 108, 109; 110, 111; 114, 115, 116; 118, 119; 137, 138, 139, 140; 141, 142; 144, 145; 146, 147; 151, 152; 153, 154.\\
\ \\
{\it 3.4.2 Lie algebras (among $1, 2,\ldots, 159)$ without $CP$'s}\\
7, 21, 22, 28, 69, 77, 100, 101; 121, 122, $\ldots$, 136; 148, 149, 150, 155 (28 cases).  This corresponds to [EO, Proposition 3.4], except for the fact that (13457H) is not a Lie algebra and that (1357S, $\xi = 1$) is isomorphic with (2357D) [G, p. 59].\\
\ \\
{\it 3.4.3. Quasi quadratic Lie algebras}\\
Namely the three quadratic ones: 7, 22, 101 (see also [FS]) and 133, 134, 135, 136 (see also [Ce2]).  All of them are coregular.\\
\ \\
{\it 3.4.4.} After consulting the list, we notice that for $\mathfrak{g}$ nonsingular the strongly complete subalgebra $M$ (and also $M_1$) of $S(\mathfrak{g})$ is already contained in $S(F(\mathfrak{g}))$.  In particular, it is also strongly complete in $S(F(\mathfrak{g}))$.  Compare this to (2) of Theorem 22.\\
\ \\
{\bf \large Acknowledgments}\\
\\
We would like to thank Michel Van den Bergh for his genuine interest and for his valuable help in calculating the invariants for some of the more difficult cases.  We are also grateful to Rudolf Rentschler and Oksana Yakimova for making some useful suggestions.\\
Special thanks go to our colleague Peter De Maesschalck for writing some indispensable and very efficient programs in MAPLE.\\
The main results of this paper (together with the ones of [O5]) were presented at the Universities of  Paris, St-Etienne, Erlangen, Poitiers and Yale University.

\section*{References}

\begin{itemize}
\bibitem[AG]{AG} J.M. Ancochea-Bermudez, M. Goze, Classification des alg\`ebres de Lie nilpotentes complexes de dimension 7, Arch. Math. 52 (1989), 157--185.
\bibitem[Bo1]{Bo1} A.V. Bolsinov, Commutative families of functions related to consistent Poisson brackets, Acta Appl. Math. 24 (1991), no 3, 253--274.
\bibitem[Bo2]{Bo2} A.V. Bolsinov, Complete commutative subalgebras in polynomial Poisson algebras: a proof of the Mishchenko-Fomenko conjecture (2008), paper available in PDF format, http://www-staff.lboro.ac.uk/~maab2/publications.html
\bibitem[BGR]{BGR} W. Borho, P. Gabriel, R. Rentschler, Primideale in Einh\"{u}llenden aufl\"{o}sbarer Lie-Algebren, Lecture Notes in Math., vol. 357, Springer-Verlag, Berlin, 1973.
\bibitem[Bou]{Bou} N. Bourbaki, Groupes et alg\`ebres de Lie, chap. I  (2nd ed.), Act. Sci. Ind. 1285, Hermann, Paris, 1971.
\bibitem[Ca]{Ca} R. Carles, Weight systems for complex nilpotent Lie algebras and application to the varieties of Lie algebras, Pr\'epublication no 96, D\'epartement de Math\'ematiques, Universit\'e de Poitiers (1996).
\bibitem[Ce1]{Ce1} A. Cerezo, Les alg\`ebres de Lie nilpotentes, r\'eelles et complexes de dimension 6, Pr\'epublication 27, D\'epartement de Math\'ematiques, Universit\'e de Nice (1983), 34p. http://math.unice.fr/~frou/ACpublications.html
\bibitem[Ce2]{Ce2} A. Cerezo, On the rational invariants of a Lie algebra, Pr\'epublication 68, D\'epartement de Math\'ematiques, Universit\'e de Nice (1985), 54p.\\ http://math.unice.fr/~frou/ACpublications.html
\bibitem[Ce3]{Ce3} A. Cerezo, Calcul des invariants alg\'ebriques et rationnels d'une matrice nilpotente,  Pr\'epublication 35, D\'epartement de Math\'ematiques, Universit\'e de Poitiers (1988), 25p. http://math.unice.fr/~frou/ACpublications.html
\bibitem[Ce4]{Ce4} A. Cerezo, Sur les invariants alg\'ebriques du groupe engendr\'e par une matrice nilpotente, unpublished, handwritten manuscript, 83 p.\\ 
http://math.unice.fr/~frou/ACinvariants/inv0.PDF   (References and Introduction, 330 Kb)\\
http://math.unice.fr/~frou/ACinvariants/inv1.PDF  (Chapter I, 1291 Kb)\\
http://math.unice.fr/~frou/ACinvariants/inv2.PDF  (Chapter II, 2567 Kb)
\bibitem[DNO]{ADNO} L. Delvaux, E. Nauwelaerts, A.I. Ooms, On the semicenter of a universal enveloping algebra, J. Algebra 94 (1985), 324--346.
\bibitem[DNOW]{DNOW} L. Delvaux, E. Nauwelaerts, A.I. Ooms, P. Wauters, Primitive localization of an enveloping algebra, J. Algebra 130 (1990), 311--327.
\bibitem[D1]{D1} J. Dixmier, Sur les repr\'esentations unitaires des groupes de Lie nilpotents, II, Bull. Soc. Math. France 85 (1957), 325--388.
\bibitem[D2]{D2} J. Dixmier, Sur les repr\'esentations unitaires des groupes de Lie nilpotents, III, Canadian J. Math., 10, (1958), 321--348.
\bibitem[D3]{D3} J. Dixmier, Sur les repr\'esentations unitaires des groupes de Lie nilpotents, IV, Canadian J. Math., 11, (1959), 321--344.
\bibitem[D4]{D4} J. Dixmier, Sur les alg\`ebres enveloppantes de $sl(n,C)$ et $af(n,C)$, Bull. Sci. Math. 100
(1976), 57--95.
\bibitem[D5]{D5} J. Dixmier, Enveloping Algebras, Grad. Stud. Math., vol 11, Amer. Math. Soc., Providence, RI, 1996.
\bibitem[DDV]{DDV} J. Dixmier, M. Duflo, M. Vergne, Sur la repr\'esentation coadjointe d'une alg\`ebre de Lie, Compositio Math. 29 (1974), 309--323.
\bibitem[E]{E} A. G. Elashvili, Frobenius Lie algebras, Funct. Anal. i Prilozhen 16 (1982) 94--95.
\bibitem[EO]{EO} A. G. Elashvili, A. I. Ooms, On commutative polarizations,  J. Algebra 264  (2003), 129--154.
\bibitem[FJ1]{FJ1} F. Fauquant-Millet, A. Joseph, Semi-centre de l'alg\`ebre enveloppante d'une sous-alg\`ebre parabolique d'une alg\`ebre de Lie semi-simple, Ann. Sci. \'Ecole Norm. Sup. 38 (2005), 155--191.
\bibitem[FJ2]{FJ2} F. Fauquant-Millet, A. Joseph, La somme des faux degr\'es - un myst\`ere en th\'eorie des invariants, Adv. Math. 217 (2008), 1476--1520.
\bibitem[FS]{FS} G. Favre, L.J. Santharoubane, Symmetric, invariant, non-degenerate bilinear form on a Lie algebra, J. Algebra 105 (1987), 451--464.
\bibitem[GK]{GK} I.M. Gelfand, A.A. Kirillov, Sur les corps li\'es aux alg\`ebres enveloppantes des alg\`ebres de Lie, Inst. Hautes Etudes Sci. Publ. Math. 31 (1966) 5--19.
\bibitem[G]{G} M-P. Gong, Classification of nilpotent Lie algebras of dimension 7, PhD Thesis, University of Waterloo, Ontario, Canada, 1998.
\bibitem[GoKh]{GoKh} M. Goze, Y. Khakimdjanov, Nilpotent Lie algebras, Mathematics and its applications, vol. 361, Kluwer, 1996.
\bibitem[GY]{GY} J. H. Grace, A. Young, The algebra of invariants, Cambridge University Press, 1903.
\bibitem[GPS]{GPS} G.-Greuel,G. Pfister, H. Sch\"{o}nemann, SINGULAR 3.0, A Computer Algebra System for Polynomial Computations, Centre for Computer Algebra, University of Kaiserslautern, 2005, http://www.singular.uni-kl.de.
\bibitem[J1]{J1} A. Joseph, Proof of the Gelfand-Kirillov conjecture for solvable Lie algebras, Proc. Amer. Math. Soc. 45 (1974) 1-–10.
\bibitem[J2]{J2} A. Joseph, A preparation theorem for the prime spectrum of a semisimple Lie algebra, J. Algebra 48 (1977) 241--289.
\bibitem[J3]{J3} A. Joseph, Parabolic actions in type A and their eigenslices, Transform. Groups, 12 (2007), no 3, 515--547.
\bibitem[J4]{J4} A. Joseph, Slices for biparabolic coadjoint actions in type A, J. Algebra 319 (2008), 5060--5100.
\bibitem[J5]{J5} A. Joseph, Compatible adapted pairs and a common slice theorem for some centralizers, Transform. Groups, 13 (2008), 637--669.
\bibitem[J6]{J6} A. Joseph, Invariants and slices for reductive and biparabolic coadjoint actions, Lecture Notes Weizmann, 20/2/2007, revised 7/1/2010,\\ http://www.wisdom.weizmann.ac.il/~gorelik/agrt.htm.
\bibitem[JL]{JL} A. Joseph, P. Lamprou, Maximal Poisson Commutative subalgebras for truncated parabolic subalgebras of maximal index in $sl_n$, Transform. Groups, 12 (2007), 549--571.
\bibitem[JS]{JS} A. Joseph, D. Shafrir, Polynomiality of invariants, unimodularity and adapted pairs. Transform. Groups, 15 (2010), 851--882.
\bibitem[K]{K} A.A. Korotkevich, Integrable Hamiltonian systems on low-dimensional Lie algebras, Sb. Math. 200, No 12, 1731--1766 (2009).
\bibitem[Ma1]{Ma1} L. Magnin, Sur les alg\`ebres de Lie nilpotentes de dimension $\leq$ 7, J. Geom. Phys., 3 (1986), 119--144.
\bibitem[Ma2]{Ma2} L. Magnin, Adjoint and trivial cohomology tables for indecomposable nilpotent Lie algebras of dimension $\leq$ 7 over C, Electronic Book, Institut de Math\'ematique, Universit\'e de Bourgogne, Second Corrected Edition (2007), (810 pages + vi), http://www.u-bourgogne.fr/monge/l.magnin
\bibitem[Ma3]{Ma3} L. Magnin, Determination of 7-dimensional indecomposable nilpotent complex Lie algebras by adjoining a derivation to 6-dimensional Lie algebras, Algebr. Represent. Theory, 13 (2010), 723--753.
\bibitem[MF]{MF} A.S. Mishchenko, A.T. Fomenko, Euler equations on finite-dimensional Lie groups, Math. USSR-Izv. 12 (1978), 371--389.
\bibitem[MW]{MW} C.C. Moore, J.A. Wolf, Square integrable representations of nilpotent groups, Trans. Amer. Math. Soc. 185 (1973), 445--462.
\bibitem[M]{M} V. Morozov,  Classification of nilpotent Lie algebras of sixth order, Izv. Vyssh. Uchebn. Zaved. Mat., 4(5), (1958), 161--171.
\bibitem[N]{N} Nghi\^em Xu\^an Hai, Sur certains sous-corps commutatifs du corps enveloppant d'une alg\`ebre de Lie r\'esoluble, Bull. Sci. Math. 96 (1972), 111--128.
\bibitem[O1]{O1} A.I. Ooms, On Lie algebras with primitive envelopes, supplements, Proc. Amer. Math. Soc. 58 (1976), 67--72.
\bibitem[O2]{O2} A.I. Ooms, On Frobenius Lie algebras, Comm. Algebra 8 (1980) 13--52.
\bibitem[O3]{O3} A.I. Ooms, On certain maximal subfields in the quotient division ring of an enveloping algebra, J. Algebra, 230 (2000), 694--712.
\bibitem[O4]{O4} A. I. Ooms, The Frobenius semiradical of a Lie algebra,  J. Algebra 273 (2004), 274--287.
\bibitem[O5]{O5} A. I. Ooms, Computing invariants and semi-invariants by means of Frobenius Lie algebras,  J. Algebra 321 (2009), 1293--1312.
\bibitem[OV]{OV} A.I. Ooms, M. Van den Bergh, A degree inequality for Lie algebras with a regular Poisson semicenter, J. Algebra 323 (2010) 305--322, arXiv: math.RT /0805.1342v1, 2008.
\bibitem[P]{P} D.I. Panyushev, On the coadjoint representation of $Z_2$-contractions of reductive Lie algebras, Adv. Math. 213 (2007), no 1, 380–-404.
\bibitem[PPY]{PPY} D.I. Panyushev, A. Premet, O. S. Yakimova, On symmetric invariants of centralizers in reductive Lie algebras, J. Algebra 313 (2007), no 1, 343--391.
\bibitem[PY]{PY} D.I. Panyushev, O. S. Yakimova, The argument shift method and maximal commutative subalgebras of Poisson algebras, Math. Res. Lett.15 (2008), no 2, 239--249. arXiv:math.RT/0702583v1 (2007).
\bibitem[PSW]{PSW} J. Patera, R. T. Sharp, P. Winternitz, H. Zassenhaus, Invariants of real low dimension Lie algebras, J. Math. Phys., 17  (1976), 986--994.
\bibitem[RV]{RV} R. Rentschler, M. Vergne, Sur le semi-centre du corps enveloppant d'une alg\`ebre de Lie, Ann. Sci. \'Ecole Norm. Sup. 6 (1973) 389--405.
\bibitem[R]{R} M. Romdhani, Classification of real and complex nilpotent Lie algebras of dimension 7, Linear and Multilinear Algebra, 24 (1989) 167--189.
\bibitem[Sa]{Sa} S. T. Sadetov, A proof of the Mishchenko-Fomenko conjecture (1981), Dokl. Akad. Nauk 397 (6) (2004), 751--754.
\bibitem[Se]{Se} C. Seeley, 7-dimensional nilpotent Lie algebras, Trans. Amer. Math. Soc., 335 (1993) 479--496.
\bibitem[Sk]{Sk} S. Skryabin, Invariants of finite group schemes, J. London Math. Soc. (2) 65 (2002), no. 2, 339--360.
\bibitem[T]{T} A.A. Tarasov, The maximality of certain commutative subalgebras in Poisson algebras of a semisimple Lie algebra, Russian Math. Surveys 57 (2002), no. 5, 1013--1014.
\bibitem[TY]{TY} P. Tauvel, R. T. W. Yu, Lie Algebras and Algebraic Groups, Springer Monographs in Mathematics, Springer-Verlag, Berlin, 2005.
\bibitem[V]{V} M. Vergne, La structure de Poisson sur l'alg\`ebre sym\'etrique d'une alg\`ebre de Lie nilpotente, Bull. Soc. Math. France, 100 (1972), 301--335.
\bibitem[Y]{Y} O. Yakimova, A counterexample to Premet's and Joseph's conjectures, Bull. London Math. Soc., 39 (2007) 749--754.
\end{itemize}
\end{document}